# RANDOM RECURRENCE EQUATIONS AND RUIN IN A MARKOV-DEPENDENT STOCHASTIC ECONOMIC ENVIRONMENT[1]

### By Jeffrey F. Collamore

*University of Copenhagen*

We develop sharp large deviation asymptotics for the probability of ruin in a Markov-dependent stochastic economic environment and study the extremes for some related Markovian processes which arise in financial and insurance mathematics, related to perpetuities and the ARCH(1) and GARCH(1, 1) time series models. Our results build upon work of Goldie [*Ann. Appl. Probab.* **1** (1991) 126–166], who has developed tail asymptotics applicable for independent sequences of random variables subject to a random recurrence equation. In contrast, we adopt a general approach based on the theory of Harris recurrent Markov chains and the associated theory of nonnegative operators, and meanwhile develop certain recurrence properties for these operators under a nonstandard "Gärtner–Ellis" assumption on the driving process.

**1. Introduction and summary.** In a variety of problems in insurance mathematics and risk management, as well as other applied areas, it is relevant to study the tail probability of a random variable satisfying a stochastic recurrence equation. An example of this type arises in risk theory, where the objective is to characterize the probability of ruin of an insurance company whose losses are governed by Lundberg's (1903) classical model, but where the company earns stochastic interest on its capital. In the setting of stochastic investments, the analysis of ruin departs substantially from the renewal-theoretic approach typically employed for Lundberg's original model. Instead, one introduces an associated process $\{W_n\}$, defined below,

Received May 2007; revised November 2008.

[1]Supported in part by Danish Research Council (SNF) Grants 21-04-0400 and 272-06-0442.

and observes that

$$\mathbb{P}\{\text{ruin}\} = \mathbb{P}\Big\{\sup_n W_n > u\Big\},$$

where $u$ is the initial capital of the company and $W := \sup_n W_n$ satisfies the *random recurrence equation*

$$(1.1) \qquad W \stackrel{d}{=} B + A \max\{0, W\}$$

for certain random variables $A$ (associated with the investment process) and $B$ (associated with the insurance business); see Section 2.1 below. The characterization of ruin then centers around (1.1) and, in particular, the tail decay of its solution as $u \to \infty$. In this setting, it is known that the probability of ruin decays at a certain polynomial rate, namely,

$$(1.2) \qquad \mathbb{P}\{W > u\} \sim Cu^{-\mathfrak{r}} \qquad \text{as } u \to \infty$$

for constants $C$ and $\mathfrak{r}$; see, for example, Goldie (1991), Nyrhinen (2001), and, for the continuous time case, Kalashnikov and Norberg (2002), Paulsen (2002), Pergamenshchikov and Zeitouny (2006) and Klüppelberg and Kostadinova (2008). Rough large deviation asymptotics in a general setting have also been developed in Nyrhinen (1999).

Related recurrence equations arise, for example, in life insurance mathematics, where attention is focused on perpetuities, which describe the discounted future payments of a life insurance company; and in financial time series modeling, where it is relevant to describe the tail decay for the now-standard ARCH(1) and GARCH(1, 1) models, used to quantify the logarithmic returns on an investment [cf. Engle (1982), Bollerslev (1986), Embrechts, Klüppelberg and Mikosch (1997)]. In these cases, the solution is obtained by solving a random recurrence equation closely related to (1.1), namely,

$$(1.3) \qquad V \stackrel{d}{=} B + AV.$$

In the case of perpetuities, the random variables $A$ and $B$ are once again determined by the investment and insurance processes, respectively, and $V$ describes the future financial obligations of the company; see Section 2.3 below. It is known—both for the case of perpetuities and for the ARCH(1) and GARCH(1, 1) financial models—that

$$(1.4) \qquad \mathbb{P}\{V > u\} \sim \tilde{C}u^{-\mathfrak{r}} \qquad \text{as } u \to \infty$$

for constants $\tilde{C}$ and $\mathfrak{r}$ [cf. Goldie (1991), Mikosch (2003), and, for a continous-time version, Carmona, Petit and Yor (2001)]. Other relevant results for perpetuities can be found, for example, in Dufresne (1990) and Cairns (1995).

These diverse problems are unified through (1.1) and (1.3), which state that the random variable of interest statisfies an equation of the general form

$$(1.5) \qquad\qquad Z \stackrel{d}{=} \Phi(Z)$$

for some real-valued random function $\Phi$. Solutions to such recurrence equations have been developed in Kesten (1973), Grincevicius (1975) and Grey (1994), and particularly Goldie (1991), who introduced an approach based on "implicit" renewal theory, which is widely applicable in the setting of (1.5). Based on the results of Goldie's paper, one readily obtains estimates such as (1.2) and (1.4), as well as various estimates relevant, for example, in queuing theory and other applied areas.

If, however, the financial or insurance process arising above is *Markov* dependent, then the above approach breaks down and the situation is actually quite different, and it is not possible to develop corresponding estimates to (1.2) and (1.4) based on the recurrence equations (1.1) and (1.3). Such extensions are nonetheless of considerable applied relevance. For example, the investment returns of an insurance company would not generally be independent, nor could they be described as a Markov chain in finite state space. A realistic model would typically involve a state space which is, say, real-valued and therefore *uncountable*, and would have increments which are *unbounded*. This is the situation in many of the standard financial models, such as the ARMA, stochastic volatility and GARCH time series models.

The objective of this article is to study (1.2) and (1.4) in a Markovian setting where the insurance and, especially, financial processes involved are driven by a Harris recurrent Markov chain. Based on the regeneration technique of Athreya and Ney (1978) and Nummelin (1978), we shall show that recurrence relations similar to (1.1) and (1.3) can be obtained. [It should be noted that this approach differs markedly from known methods for similar problems; cf., e.g., Nyrhinen (2001) and de Saporta (2005).] In a general Markovian setting, the analysis of these equations turns out, however, to be considerably more complicated than in the independent case, and a main aspect of our study will be centered upon the regularity properties of certain random quantities formed over the regeneration cycles of the Markov chain. In Theorem 4.2 below—a central result of this paper—we develop regularity properties closely related to geometric $\mathfrak{r}$-recurrence for the operator $\hat{P}_{\mathfrak{r}}$, where

$$\hat{P}_{\alpha}(x, dy) := e^{\alpha f(y)} P(x, dy) \qquad \forall \alpha,$$

and $P$ is the transition kernel of the underlying Markov chain and $f$ a real-valued function. Geometric $\mathfrak{r}$-recurrence plays an important role in the large deviations theory for general Markov chains [cf. Ney and Nummelin (1987a, 1987b)].

In particular, there has been much recent attention focused on establishing geometric recurrence for a given Markov transition kernel $P$, which is typically achieved by verifying a Lyapunov drift condition, namely,

$$(\mathfrak{D}) \qquad \int V(y)P(x,dy) \le \rho V(x) + b\mathbf{1}_C(x)$$

for some function $V \ge 1$, "small set" $C$ and constants $b < \infty$ and $\rho < 1$ [cf. Meyn and Tweedie (1993), Chapter 15]. However, validating a condition such as $(\mathfrak{D})$ with the *operator* $\hat{P}_\tau$ in place of $P$—which would yield geometric $\tau$-recurrence and some extensions developed here—is generally much more difficult. In this paper, we propose an alternative approach, based on the introduction of an auxiliary "$h$-function" and the verification of an associated "Gärtner–Ellis" limit and minorization condition. These conditions provide an alternative to $(\mathfrak{D})$ which, at least in the context of our examples, can be verified somewhat more naturally. [An alternative approach has recently been introduced in Kontoyiannis and Meyn (2003, 2005). The relationship between their approach and ours is discussed in more detail in Section 2 below.] In our development, we shall rely heavily on the theory of nonnegative operators, as summarized in Nummelin (1984).

We conclude by mentioning some recent work on Markov-driven perpetuities, as described in (1.4) above. An extension of (1.4) to finite state space Markov chains has recently been obtained in de Saporta (2005). Also, an extension to functionals of continuous-time Markov processes has been given in Blanchet and Glynn (2005), but under a boundedness assumption on the functionals which is violated in the examples we consider here. Specifically, their assumptions in our problem would imply that the sequence $\{\sum_{i=1}^n A_i\}$ is bounded, leading to exponential rather than polynomial decay for the probability of ruin. In fact, our emphasis on unbounded processes $\{A_n\}$ will lead to the main technical difficulties which we shall encounter below. In contrast to both of these papers, our methods will be based on the regeneration properties of the underlying Markov chain, and regeneration will be central to our approach here.

An outline of this paper is as follows. In the next section, we give a more precise description of the ruin problem with stochastic investments, as introduced in (1.1), (1.2), and then turn to a Markovian formulation of this problem. Next, we describe the same Markovian formulation, but for perpetuities and the ARCH(1) and GARCH(1,1) financial models. Some examples are given in Section 3, proofs are given in Sections 4 and 5 and generalizations are briefly discussed in Section 6.

## 2. Statement of results.

2.1. *A description of the insurance risk model in the i.i.d. setting.* We begin by recalling the classical Cramér–Lundberg model for the capital growth of an insurance company, namely,

$$Y_t = u + ct - \sum_{i=1}^{N_t} \xi_i, \tag{2.1}$$

where $u$ is the initial capital of the company, $c$ is the premiums income, $\{\xi_i\}$ are the i.i.d. claims losses and $\{N_t\}$ is a Poisson process, independent of $\{\xi_i\}$, which describes the occurrence times of the claims. These assumptions imply that $\{Y_t\}$ is a Lévy process, with i.i.d. losses over unit intervals given by

$$B_n := -(Y_n - Y_{n-1}), \qquad n = 1, 2, \dots.$$

It is assumed that $\{Y_t\}_{t \geq 0}$ has a positive drift or, equivalently, that $\{B_n\}_{n \in \mathbb{Z}_+}$ has a negative mean.

We now depart from this model by introducing a financial process describing the investment returns. Assume that the return rate during the $n$th discrete time interval is given by $r_n$, and let $Z_n$ denote the total capital of the insurance company at time $n$. Then

$$Z_n = (1 + r_n)(Z_{n-1} - B_n), \qquad n = 1, 2, \dots \quad \text{and} \quad Z_0 = u. \tag{2.2}$$

[One could alternatively assume that $Z_n = (1 + r_n)Z_{n-1} - B_n$, and the analysis would carry through without significant change.] Setting $R_n = 1 + r_n$ and solving (2.2) recursively for $Z_n$ yields

$$\begin{aligned} Z_n &= R_n(Z_{n-1} - B_n) = R_n R_{n-1}(Z_{n-2} - B_{n-1}) - R_n B_n \\ &= \dots = (R_n \cdots R_1)Z_0 - (R_n \cdots R_1)B_1 - \dots - R_n B_n. \end{aligned} \tag{2.3}$$

Assuming $r_n > -1$ a.s. for all $n$, we may set $A_n := 1/(1 + r_n)$ to be the stochastic discount factor, and multiplying left- and right-hand sides of (2.3) by $A_1 \cdots A_n$ yields

$$(A_1 \cdots A_n)Z_n = u - W_n, \tag{2.4}$$

where

$$W_n := B_1 + A_1 B_2 + \dots + (A_1 \cdots A_{n-1})B_n, \tag{2.5}$$

that is, $W_n$ represents the total discounted losses incurred at times $1 \leq i \leq n$.

Let

$$W = \sup_{n \geq 1} W_n.$$

Then $\{Z_n < 0, \text{ for some } n\} \iff \{W > u\}$, by (2.4), and hence the probability of ruin is given by

$$\Psi(u) := \mathbb{P}\{W > u\}. \tag{2.6}$$

2.2. *The Markovian case and a statement of our results.* Our next objective is to generalize the above formulation so that it allows for Markov dependence. In the interest of simplicity we shall first assume that there is dependence only in the investment process $\{A_n\}$, and that the insurance process $\{B_n\}$ is i.i.d. and independent of $\{A_n\}$. In Section 6.1 below, we shall discuss a slight generalization which allows for Markov dependence also in the sequence $\{B_n\}$.

We begin, then, with a Markov chain in general state space, denoted by $\{X_n\}$, and assume

$$\log A_n = f(X_n), \tag{2.7}$$

where $f : \mathbb{S} \to \mathbb{R}$ and is typically *unbounded.* (We could equally well assume that the function $f$ is random; see Remark 2.3 below.)

We suppose that $\{X_n\}$ is time-homogeneous, taking values in a countably generated measurable state space $(\mathbb{S}, \mathcal{S})$ with transition kernel $P(x, E)$ and $k$-step transition kernel $P^k := PP^{k-1}$ for $k > 1$. Assume that $\{X_n\}$ is aperiodic and irreducible with respect to a maximal irreducibility measure $\varphi$. [For the definitions and a further characterization of these conditions, see Nummelin (1984) or Meyn and Tweedie (1993).] As an additional regularity condition, we suppose that for $\varphi$-a.a. initial states $x \in \mathbb{S}$, the distribution of $\sum_{i=1}^{n} \log A_i$ is spread out for all $n \geq N(x)$, where $N(x)$ is some positive integer.

Under the assumption that $\{X_n\}$ is irreducible with respect to $\varphi$, there exists a *minorization* for $\{X_n\}$ [Nummelin (1984), Theorem 2.1]; namely,

$$(\mathcal{M}_0) \qquad \delta \mathbf{1}_C(x) \nu(E) \leq P^k(x, E) \qquad \forall x \in \mathbb{S}, \ E \in \mathcal{S},$$

for some $k \in \mathbb{Z}_+$, $\delta \in (0, 1]$, a probability measure $\nu$ on $(\mathbb{S}, \mathcal{S})$ and a "small set" $C \in \mathcal{S}$, where $\varphi(C) > 0$. Here we will work with a strengthening of this property, condition $(\mathfrak{M})$, which will be described below.

*Notation.*

$$S_n = \log A_1 + \cdots + \log A_n, \qquad n = 1, 2, \ldots \quad \text{and} \quad S_0 = 0;$$

$$S_n^{(h)} = h(X_1) + \cdots + h(X_n), \qquad n = 1, 2, \ldots,$$

$$\text{for any measurable function } h : \mathbb{S} \to [0, \infty);$$

$$\Lambda(\alpha) = \limsup_{n \to \infty} \frac{1}{n} \log \mathbb{E}[e^{\alpha S_n}] \qquad \forall \alpha \in \mathbb{R};$$

$$\Lambda_B(\alpha) = \log \mathbb{E}[|B_1|^\alpha] \qquad \forall \alpha \in \mathbb{R};$$

$$\boldsymbol{\Lambda}(\boldsymbol{\alpha}) = \limsup_{n \to \infty} \frac{1}{n} \log \mathbb{E}[e^{\alpha S_n + \beta S_n^{(h)}}] \qquad \forall \boldsymbol{\alpha} = (\alpha, \beta) \in \mathbb{R}^2;$$

$$\mathcal{L}_a h = \{x \in \mathbb{S} : h(x) \leq a\} \qquad \forall a \in \mathbb{R}, \text{ for any function } h : \mathbb{S} \to \mathbb{R}.$$

Also let $\mathbf{1}_C$ denote the indicator function on $C$, and for any measure $\nu$, let supp $\nu$ denote the support of $\nu$. (In the notation for $\Lambda$ and $\boldsymbol{\Lambda}$, we have suppressed the dependence on the initial state of the Markov chain. However, in Proposition 5.1 below, it will be shown that these quantities are actually independent of this initial state.)

The functions $\Lambda$ and $\boldsymbol{\Lambda}$ are the "Gärtner–Ellis" limits arising in large deviation theory [cf. Dembo and Zeitouni (1998), Chapter 2]. Roughly, $\Lambda$ can be equated to the spectral radius of the transform kernel

$$(2.8) \qquad \hat{P}_\alpha(x, E) := \int_E e^{\alpha f(y)} P(x, dy) \qquad \forall x \in \mathbb{S}, \ E \in \mathcal{S},$$

where $f$ is given as in (2.7), and similarly for $\boldsymbol{\Lambda}$; cf. de Acosta (1988), Section 7.

We turn now to some assumptions on the Markov-additive process $\{(X_n, S_n) : n = 0, 1, \ldots\}$. First note that, if the average return rate is *positive*, then $1 + r_1 > 1$ "on average," and hence it is reasonable to expect that

$$\mathbb{E}_\pi[A_1] := \mathbb{E}_\pi[(1 + r_1)^{-1}] \in (0, 1),$$

where $\pi$ is the stationary probability measure of $\{X_n\}$. By Ney and Nummelin (1987a), Lemmas 3.3 and 5.2 (and their proofs) and an application of Jensen's inequality, we then obtain under very weak regularity conditions on the Markov chain that, if $\Lambda$ is finite in a neighborhood of zero,

$$(2.9) \qquad \Lambda'(0) = \mathbb{E}_\pi[\log A_1] < 0.$$

Hence

$$(2.10) \qquad \mathfrak{r} := \sup\{\alpha : \Lambda(\alpha) \le 0\} > 0.$$

Moreover, if $\mathbb{P}_\pi\{A_1 > 1\} > 0$, then $\mathfrak{r} < \infty$, so the solution to (2.10) is in fact a proper solution.

In addition to the existence of the solution in (2.10), we will also need to assume a further regularity condition on $\Lambda$, roughly stating that it is finite in a neighborhood of $\mathfrak{r}$, and that for some choice of $h$, the generating function of $\{h(X_n)\}$ is sufficiently well behaved around the origin. To motivate this condition, set

$$h(x) = |f(x)| \quad \text{and} \quad S_n^- = \sum_{i=1}^n (-f(X_i)) \vee 0$$

and observe by Hölder's inequality that

$$(2.11) \qquad \mathbb{E}[e^{\alpha S_n + \beta S_n^{(h)}}] \le (\mathbb{E}[e^{p(\alpha + \beta) S_n}])^{1/p} (\mathbb{E}[e^{2q\beta S_n^-}])^{1/q},$$

where $p^{-1} + q^{-1} = 1$. Then

$$(2.12) \qquad \boldsymbol{\Lambda}((\alpha, \beta)) \le \frac{1}{p} \Lambda(\tilde{\alpha}) + \frac{1}{q} \Lambda(-\varepsilon),$$

where $\tilde{\alpha} = p(\alpha + \beta)$ and $\varepsilon = 2q\beta$. Choosing $(\alpha, \beta)$ sufficiently close to $(\mathfrak{r}, 0)$, then $p$ sufficiently close to one, and finally letting $\beta \searrow 0$, we conclude that if $\Lambda$ is finite in a neighborhood of the interval $[0, \mathfrak{r}]$, then

$$(2.13) \qquad \mathbf{\Lambda}((\alpha, \beta)) < \infty \qquad \text{for some } \alpha > \mathfrak{r} \text{ and } \beta > 0.$$

More generally, if we take $h$ to be an arbitrary function, then in place of (2.11) we obtain

$$(2.14) \qquad \mathbb{E}[e^{\alpha S_n + \beta S_n^{(h)}}] \le (\mathbb{E}[e^{p\alpha S_n}])^{1/p} (\mathbb{E}[e^{q\beta S_n^{(h)}}])^{1/q}$$

and hence (2.13) still holds, provided that $\Lambda$ is finite in a neighborhood of $\mathfrak{r}$ and the generating function of $\{S_n^{(h)}\}$ is finite in a neighborhood of zero. Hence, we may also choose $h(x) = \|x\|$ or a more slowly increasing function, such as $\log \|x\| \vee 0$. Later, we will relate the function $h$ to the minorization condition $(\mathfrak{M})$, given below, and for this reason it will often be necessary to choose $h$ to be *different* from $|f|$. This is because we will need $h(x)$ to tend to infinity as $\|x\| \nearrow \infty$ in order for $(\mathfrak{M})$ to be satisfied. See, for example, Example 3.3 below. (Also see Example 3.1 for another case where we would generally not choose $h = |f|$.)

The above considerations provide motivation for the following.

*Hypotheses.*

(H$_1$)  $\mathfrak{r} \in (0, \infty)$.

(H$_2$)  There exists a function $h \colon \mathbb{S} \to \mathbb{R}$ and points $\alpha > \mathfrak{r}$ and $\beta > 0$ such that $\mathbf{\Lambda}((\alpha, \beta)) < \infty$ and $\Lambda_B(\alpha) < \infty$.

Next we introduce a further condition which regulates the behavior of $\{A_n\}$ and $\{X_n\}$ on the *level sets* of the function $h$ appearing in (H$_2$). To motivate this condition, note that for a very large class of financial processes, there is some form of stochastic domination. A simple example of this type is the AR(1) process,

$$(2.15) \qquad X_n = cX_{n-1} + \zeta_n, \qquad n = 1, 2, \ldots \quad \text{and} \quad X_0 = x,$$

for $|c| \in (0, 1)$ and $\{\zeta_n\}$ an i.i.d. sequence of standard Gaussian random variables. Then

$$X_0 = z \quad \Longleftrightarrow \quad X_1 \sim \text{Normal}(cz, 1).$$

Consequently, if $c > 0$ then

$$(2.16) \quad x \le y \quad \Longrightarrow \quad P(x, E) \le P(y, E) \qquad \text{for all sets } E \subseteq \left[\frac{c}{2}(x + y), \infty\right);$$

and conversely,

$$x \leq y \quad \implies \quad P(y, E) \leq P(x, E)$$

(2.17)

$$\text{for all sets } E \subseteq \left(-\infty, \frac{c}{2}(x+y)\right],$$

where $P$ denotes the Markov transition kernel of $\{X_n\}$. Therefore, if $\mathcal{L}_a h \subseteq [-b, b]$ then

(2.18) $\quad P(x, E) \leq P(-b, E) + P(b, E) \qquad \text{for all } x \in \mathcal{L}_a h \text{ and } E \in \mathcal{S}.$

If $c < 0$, then the inequalities on the right-hand sides of (2.16) and (2.17) are reversed, but (2.18) remains valid. This type of reasoning can be generalized to include, for example, general ARMA$(p, q)$ models by using a vector representation for the process $\{X_n\}$; cf. Meyn and Tweedie (1993), Chapter 2 and Section 6 below. More complicated financial processes can be handled similarly.

Of course, a representation such as (2.18) would be quite meaningless if it were only to hold for the individual points $b$ and $-b$, which have $\varphi$-measure zero, whereas (2.18) actually holds for every $\tilde{b} \geq b$. For this reason, in (2.18) it is natural to substitute *sets* of positive $\varphi$-measure for the individual points $-b$ and $b$; and in this way we arrive at the following general condition.

(H$_3$) For any $a > 0$, there exist $\varphi$-positive sets $E_1, \ldots, E_l \subseteq \mathbb{S}$, possibly dependent on $a$, and a finite constant $D_a$ such that

$$P(x, E) \leq D_a \inf\left\{\sum_{i=1}^{l} P(x_i, E) : x_i \in E_i, 1 \leq i \leq l\right\} \qquad \forall x \in \mathcal{L}_a h, \ E \in \mathcal{S}.$$

Finally, we introduce a strengthening of the minorization ($\mathcal{M}_0$) described above.

MINORIZATION ($\mathfrak{M}$). For any $a > 0$ sufficiently large, there exists a constant $\delta_a > 0$ and a probability measure $\nu_a$ on $(\mathbb{S}, \mathcal{S})$ with $\nu_a(\mathcal{L}_a h) > 0$ such that

($\mathcal{M}$) $\qquad \delta_a \mathbf{1}_{\mathcal{L}_a h}(x)\nu_a(E) \leq P(x, E) \qquad \forall x \in \mathbb{S}, E \in \mathcal{S}.$

It should be emphasized that the function $h$ in ($\mathcal{M}$) is the same function as that appearing in (H$_2$). Thus, while (H$_2$) benefits from a small choice of $h$, the minorization ($\mathcal{M}$) benefits from a large choice of $h$ and, in practice, a balance is needed in selecting a proper choice for this function.

REMARK 2.1. There is no loss of generality in assuming that $\operatorname{supp} \nu_a \subseteq \mathcal{L}_a h \cap \mathcal{L}_b f$, where $b \in \mathbb{R}$ is arbitrarily large, since we may always truncate the measure $\nu_a$ in ($\mathcal{M}$) and the minorization will still hold. In the sequel, it will always be assumed that $\nu_a$ has been chosen in this manner.

REMARK 2.2. For simplicity, we have taken $k = 1$ in ($\mathcal{M}$) [compare ($\mathcal{M}_0$)], which may be restrictive in certain examples. A generalization to the case where $k > 1$ will be discussed below in Section 6.2.

To see how ($\mathfrak{M}$) relates to some more standard conditions, suppose for the moment that $\{X_n\}$ is uniformly recurrent, that is,

$$\delta \nu_0(E) \le P(x, E) \le d \nu_0(E)$$

($\mathfrak{R}$)

$$\forall x \in \mathbb{S}, \ E \in \mathcal{S}, \text{for some probability measure } \nu_0.$$

In this case, the $C$-set in ($\mathcal{M}_0$) may be taken to be the entire state space, which means that ($\mathfrak{M}$) then holds with $h \equiv 1$. Incidentally, from ($\mathfrak{R}$) we also obtain

$$(2.19) \qquad P(x, E) \le \frac{d}{\delta} P(x_0, E) \qquad \forall x \in \mathbb{S}, \ E \in \mathcal{S}, \text{ for any } x_0 \in \mathbb{S}.$$

Consequently (H$_3$) holds. Thus, in this setting, our conditions reduce essentially to (H$_2$) with $h \equiv 1$, namely, we require that $\Lambda(\alpha) < \infty$ and $\Lambda_B(\alpha) < \infty$, for some $\alpha > \mathfrak{r}$.

In the examples below, the upper and lower bounds in ($\mathfrak{R}$) will *not* be satisfied, but in place of the lower bound we will have

($\mathcal{M}'$) $\qquad \delta_j \mathbf{1}_{C_j}(x) \nu_j(E) \le P(x, E) \qquad \forall x \in \mathbb{S}, \ E \in \mathcal{S}, \ j \in \mathbb{Z}_+,$

along an appropriately chosen sequence of sets $C_j \nearrow \mathbb{S}$. Then an unbounded function $h$ may essentially always be found which satisfies condition ($\mathfrak{M}$), although this imposes an additional constraint on (H$_2$) as compared with the case $h \equiv 1$. For example, a typical choice for $h$ would be to take $h(x) = \|x\|$, in which case the level sets $\mathcal{L}_a h$ would tend to $\mathbb{S}$ as $a \to \infty$. In essence, then, ($\mathfrak{M}$) requires that the minorization ($\mathcal{M}_0$) hold for an arbitrarily large choice of $C$, but it is *not* required that $C = \mathbb{S}$. In particular, our assumptions lie somewhere between the minorization condition ($\mathcal{M}_0$), which is essentially always valid, and the much stronger condition that the Markov chain be uniformly recurrent. Nonetheless, these conditions are flexible enough to handle some reasonably complicated financial models, as will be seen in Section 3 below.

We note that if $P$ were replaced with $P^{k_j}$ on the right-hand side of ($\mathcal{M}'$), then this condition would indeed be an exceedingly mild requirement [cf. Meyn and Tweedie (1993), Proposition 5.2.4].

THEOREM 2.1. *Let $\{(X_n, S_n) : n = 0, 1, \dots\}$ be a Markov-additive process and let $h : \mathbb{S} \to [0, \infty)$ be a function satisfying conditions* (H$_1$)–(H$_3$) *and ($\mathfrak{M}$). Then for $\varphi$-a.a. $x \in \mathbb{S}$,*

$$(2.20) \qquad \lim_{u \to \infty} u^{\mathfrak{r}} \Psi_x(u) = C,$$

*where $C \in [0, \infty)$, and $\mathfrak{r} \in (0, \infty)$ is given as in (2.10).*

REMARK 2.3. As an extension, one could suppose that the function $f$ in (2.7) is random. In Ney and Nummelin (1987a, 1987b), Markov-additive processes are studied which take the general form $S_n = \xi_1 + \cdots + \xi_n$, where

$$\mathbb{P}\{(X_n, \xi_n) \in E \times \Gamma | \mathfrak{F}_{n-1}\} = \mathbb{P}\{(X_n, \xi_n) \in E \times \Gamma | X_{n-1}\} = \int_E P(x, dy)\mathcal{Q}(y, \Gamma)$$

for some family of probability measures $\{\mathcal{Q}(x, \Gamma) : x \in \mathbb{S}, \Gamma \in \mathcal{R}\}$, where $\mathcal{R}$ denotes the Borel $\sigma$-algebra on $\mathbb{R}$ and $\mathfrak{F}_n = \sigma\{X_0, \ldots, X_n, \xi_1, \ldots, \xi_n\}$. Thus, with a slight abuse of notation, we may write $\xi_n = f_n(X_n)$, where $\{f_n(x) : x \in \mathbb{S}, n = 1, 2, \ldots\}$ is a family of independent random variables, also independent of $\{X_n\}$, whose elements have, for fixed $x$, a common distribution function. But then $\{(X_n, \xi_n)\}$ is itself a Markov chain, which inherits a minorization from $\{X_n\}$, and clearly $\xi_n = f(X_n, \xi_n)$ for the deterministic function $f(x, y) = y$. In short, the introduction of a random function in (2.7) does not lead to additional generality, at least in principle, and the previous theorem could also have been phrased at that level of generality.

REMARK 2.4. A precise representation for the constant $C$ can be ascertained from the proof in Section 4. Under very weak conditions, it can be shown that this constant is positive. However, we will not explore the precise conditions here.

2.3. *Further remarks on our hypotheses.* Before turning to our next result we would first like to comment, briefly, on the comparison of our approach to some other methods in the literature.

An alternative approach would be to replace $(H_2)$ and $(H_3)$ with a weaker Lyapunov drift condition, namely,

$$(\mathfrak{D}_1) \qquad \int_{\mathbb{S}} V(y)e^{\mathfrak{r}f(y)}P(x, dy) \leq \rho V(x) + b\mathbf{1}_C(x),$$

where $C$ is a small set, say, and $\rho < 1$. Indeed, we shall actually utilize $(H_2)$ and $(H_3)$ together with the inherent eigenvalue and eigenfunction properties of the chain to deduce

$$(\mathfrak{D}_2) \qquad \int_{C^c} V(y)e^{\mathfrak{r}f(y)}P(x, dy) \leq \rho V(x) \qquad \forall x \in \mathbb{S},$$

where $C = \mathcal{L}_a h$ is the small set in $(\mathcal{M})$, and this is roughly equivalent to $(\mathfrak{D}_1)$. [More precisely, if $\mathfrak{r} = 0$ then either condition may be used to obtain geometric recurrence; cf. Nummelin (1984), Meyn and Tweedie (1993). It should be remarked that we will actually establish and apply a condition slightly stronger that $(\mathfrak{D}_2)$.] In this paper, the deduction of $(\mathfrak{D}_2)$ will be obtained by an indirect argument, which will constitute an important part of the proof of Theorem 4.2 below.

Nonetheless, ($\mathfrak{D}_1$) could also be viewed as a possible starting point for our results, and this approach has recently been followed, for example, by Chan and Lai (2007) and Balaji and Meyn (2000). For countable state space chains Balaji and Meyn (2000) have shown that ($\mathfrak{D}_1$) is equivalent to geometric recurrence for the $\mathfrak{r}$-shifted chain. This $\mathfrak{r}$-shifted chain will also appear in our analysis, in Section 5, and in Theorem 4.2 we shall develop similar recurrence properties, although our methods and the exact statement of our results will be quite different. Specifically, instead of geometric recurrence, we shall establish certain related moment properties. The connection between these two notions will be explained in more detail in Section 5.

In a comprehensive study, Kontoyiannis and Meyn (2003, 2005) have considered certain extensions of Balaji and Meyn (2000) to general state space chains. Specifically, they have developed multiplicative mean ergodic theory and its connection with the multiplicative Poisson equation. These results relate closely to the existence and characterization of the eigenvalues and eigenfunctions associated with the kernel $\hat{P}_\alpha$. We refer the reader to Kontoyiannis and Meyn (2003, 2005), where a survey of some other related results can also be found.

In the context of our examples, it often seems more natural to verify (H$_2$) and (H$_3$) than a condition such as ($\mathfrak{D}_1$), assuming that the driving Markov chain has finite exponential moments around the origin and therefore fits within the framework typically studied in modern large deviations theory. As mentioned in the previous section, a sufficient condition for (H$_2$) to hold is the finiteness of the Gärtner–Ellis limit, $\Lambda$, in a neighborhood of $\mathfrak{r}$, and the finiteness of an associated limit for $\{S_n^{(h)}\}$ in a neighborhood of zero, where necessarily $\Lambda(\mathfrak{r}) < \infty$ if the conclusions of our main results are to hold. While it is often difficult to obtain closed-form expressions for Gärtner–Ellis limits, their finiteness can frequently be verified by indirect means. In Example 3.3 below, finiteness is obtained along the positive axis due to the boundedness from above of the function $f$; whereas in this example, it would not be transparent that any function $V$ should satisfy ($\mathfrak{D}_1$). (This distinction would be even more striking if the interest process in that example were taken to be stochastic.) In any case, our approach exposes an interesting interplay between the Gärtner–Ellis limit of large deviations theory and geometric recurrence of the kernel $\hat{P}_\alpha$, and the latter property has important implications in the large deviations theory for Markov chains; cf. Ney and Nummelin (1987a, 1987b).

2.4. *Perpetuities and the* GARCH$(1, 1)$ *process.* A related but simpler problem to the one considered in Section 2.2 is the study of perpetuities. Assume for the moment that $\{(A_n, B_n) : n = 1, 2, \ldots\}$ is an i.i.d. sequence of random variables, and consider the tail of $\lim_{n \to \infty} W_n$ where, as before,

$$(2.21) \qquad W_n = B_1 + A_1 B_2 + \cdots + (A_1 \cdots A_{n-1}) B_n.$$

In the context of life insurance mathematics, the sequence $\{B_n\}$ typically denotes the future payments from an insurance company to its policy holders (or vice versa) at times $n = 1, 2, \ldots$, while $\{A_n\}$ denotes the discount factors associated with the investment returns, that is, $A_n = (1 + r_n)^{-1}$, where $r_n$ is the return rate at time $n$. Here, both the processes $\{A_n\}$ and $\{B_n\}$ are assumed to be random. Then $(A_1 \cdots A_{n-1})B_n$ denotes the amount of capital which needs to be set aside to cover payments at time $n$ to the policy holders, and

$$(2.22) \qquad W_\infty := \lim_{n \to \infty} W_n$$

represents the company's total future financial commitment. Under our hypotheses, the a.s. existence of the limit (2.22) follows from Goldie and Maller (2000), Theorem 2.1.

An analogous mathematical problem arises when characterizing the extremal behavior of the ARCH(1) and GARCH(1,1) financial time series models. In the GARCH(1,1) model, the logarithmic return of a stock at time $n$, denoted $R_n^*$, is governed by the system of equations

$$(2.23) \quad R_n^* = \sigma_n \xi_n, \qquad \text{where } \sigma_n^2 = a_0 + b_1 \sigma_{n-1}^2 + a_1 (R_{n-1}^*)^2, \ n = 1, 2, \ldots,$$

for $\{\xi_n\}$ an i.i.d. sequence of standard Gaussian random variables, where $a_0, a_1$ and $b_1$ are positive constants. Setting $W_n^* = \sigma_n^2$ gives

$$(2.24) \quad W_n^* = A_n W_{n-1}^* + B_n, \qquad n = 1, 2, \ldots, \quad \text{and} \quad W_0^* = y \in \mathbb{R},$$

where $A_n = b_1 + a_1 \xi_{n-1}^2$ and $B_n = a_0$. Solving (2.24) yields

$$(2.25) \quad W_n^* = (A_n \cdots A_1) W_0^* + (A_n \cdots A_2) B_1 + \cdots + A_n B_{n-1} + B_n,$$

which has a similar, although not identical, form to (2.21). Then it is of interest to study the tail of $W_\infty^* := \lim_{n \to \infty} W_n^*$. [Here we consider the limit in law. The existence of the limit distribution is then guaranteed by Goldie and Maller (2000), Theorem 3.1.] For a more detailed description of the ARCH(1) and GARCH(1,1) financial time series models, see, for example, Embrechts, Klüppelberg and Mikosch (1997) or Mikosch (2003). What the random variables $W_\infty$ and $W_\infty^*$ have in common is that they both satisfy the random recurrence equation

$$(2.26) \qquad V \stackrel{d}{=} B + AV,$$

where $(A, B) \stackrel{d}{=} (A_1, B_1)$.

It is of theoretical and applied interest to consider (2.21) and (2.25) in a setting where the process $\{A_n\}$ represents a general Markov-dependent sequence of random variables. Let

$$\tilde\Psi_x(u) = \mathbb{P}\{W_\infty > u | X_0 = x\} \quad \text{and} \quad \tilde\Psi_x^*(u) = \mathbb{P}\{W_\infty^* > u | X_0 = x\}.$$

If we adopt the same assumptions as in Section 2.2, then as a natural variant of Theorem 2.1 we obtain the following.

THEOREM 2.2.   *Let $\{(X_n, S_n) : n = 0, 1, \ldots\}$ be a Markov-additive process and let $h : \mathbb{S} \to [0, \infty)$ be a function satisfying conditions $(\mathrm{H}_1)$–$(\mathrm{H}_3)$ and $(\mathfrak{M})$. Then for $\varphi$-a.a. $x$,*

$$\lim_{u \to \infty} u^{\mathfrak{r}} \tilde{\Psi}_x(u) = \tilde{C}, \tag{2.27}$$

*where $\tilde{C} \in [0, \infty)$, and $\mathfrak{r} \in (0, \infty)$ is given as in (2.10). Moreover, (2.27) also holds if $\tilde{\Psi}_x(u)$ is replaced with $\tilde{\Psi}_x^*(u)$.*

## 3. Examples.

The objective of this section is to relate our theorems and conditions to some standard processes arising in insurance and financial mathematics.

From a mathematical point of view, it should first be noted that if the Markov chain has finite state space or is uniformly recurrent, then our conditions hold without further restrictions on the Markov chain. [However, in that setting, the proofs of our main results could be simplified considerably.] Our primary objective is to consider the case of general, Harris recurrent chains, and this setting is indeed realistic from a financial perspective, as the most reasonable models often involve a more intricate dependence structure than can be described with, say, a finite-state chain. Nonetheless, we will begin with the finite-state case and first illustrate our conditions in that context before turning to some more complicated examples.

EXAMPLE 3.1.   Assume that an insurance company receives premiums and incurs claims according to the classical Cramér–Lundberg model described in Section 2.1; thus the one-period losses, $\{B_n\}$, form an i.i.d. sequence of random variables. Suppose that the company invests its excess capital and that the investment returns are governed by the standard Black–Scholes model, but modified to allow for regime switching. More precisely, assume that the investment returns follow the stochastic differential equation

$$d\mathcal{S}_t = \sum_{i=0}^{1} \mathcal{S}_t(\mu_i \, dt + \sigma_i \, dW_t) \mathbf{1}_{\{X_{\lfloor t \rfloor} = i\}}, \tag{3.1}$$

where $\{X_n\}$ is a discrete-time Markov chain taking values in $\{0, 1\}$, $\{W_t\}$ is standard Brownian motion, $\mu_i$ and $\sigma_i$ are constants for each $i$, and it is assumed that the processes $\{W_t\}$ and $\{X_n\}$ are independent. Here the motivation is that, under different external conditions, represented by "states," the returns on the investments change. Integration of (3.1) then yields

$$A_n := \frac{\mathcal{S}_{n-1}}{\mathcal{S}_n} = \mathbf{1}_{\{X_n = 0\}} A_n^{(0)} + \mathbf{1}_{\{X_n = 1\}} A_n^{(1)}, \tag{3.2}$$

where

$$(3.3) \quad A_n^{(i)} = \exp\left\{-\left(\mu_i - \frac{\sigma_i^2}{2}\right) - \sigma_i Z_n\right\}, \qquad \{Z_n\} \sim \text{i.i.d. Normal}(0,1).$$

Regime-switching models have been introduced, for example, in Hamilton (1989), and have subsequently appeared rather widely in the literature; see, for example, Hardy (2001). (Such models could also allow for dependence between the underlying Markov chain and the conditional returns on the investments.)

A slightly more general model than (3.2) would be to assume that $\{A_n\}$ is modulated by a $k$-state Markov chain and that the discounted returns, conditional on $\{X_n\}$, are independent [but not necessarily dictated by the distribution described in (3.3)]. Namely, assume

$$(3.4) \qquad A_n = \sum_{i=0}^{k-1} \mathbf{1}_{\{X_n = i\}} A_n^{(i)},$$

where $\{(A_n^{(0)}, \dots, A_n^{(k-1)})\}$ is an i.i.d. sequence of random variables and $\{X_n\}$ is a Markov chain on $\{0, 1, \dots, k-1\}$.

In this setting, it is easy to verify our conditions based on the observation that a finite-state Markov chain satisfies a condition tantamount to uniform recurrence, namely,

$$(\mathfrak{R}) \qquad \delta\nu_0(E) \le P(x, E) \le d\nu_0(E) \qquad \forall x \in \mathbb{S}, \ E \in \mathcal{S},$$

for some probability measure $\nu_0$ and certain positive constants $\delta$ and $d$. Therefore, as already noted in the discussion following Remark 2.2, we immediately obtain $(H_3)$, and we obtain $(\mathcal{M}_0)$ with $k = 1$ and $C = \mathbb{S}$, that is, we may choose $C$ to be the *entire* state space of the Markov chain. Because $(\mathcal{M}_0)$ holds with $C = \mathbb{S}$, we consequently conclude that $(\mathfrak{M})$ holds with $h \equiv 1$. Choosing $h \equiv 1$, we see that a sufficient condition for $(H_2)$ to hold is that $\Lambda$ is finite in a neighborhood of the set $\{\alpha : \Lambda(\alpha) \le 0\}$, which is a rather weak requirement. Finally, it is well known that (2.9) holds under $(\mathfrak{R})$ [cf. Iscoe, Ney and Nummelin (1985)], and thus $(H_1)$ holds provided that $\mathbb{E}_\pi[\log A_1] < 0$, where $\pi$ is the stationary measure of the Markov chain.

The above reasoning also applies in general state space if the underlying Markov chain is uniformly recurrent, but this reasoning will *fail* in the absence of uniform recurrence. Specifically, in that case we will *not* be able to choose $h \equiv 1$.

For a discussion of uniformly recurrent chains see, for example, Iscoe, Ney and Nummelin (1985) and references therein.

EXAMPLE 3.2. Assume again the Cramér–Lundberg model for the insurance business, but now assume that the discounted logarithmic investment returns are modeled as an AR(1) process with negative drift; more precisely,

$$(3.5) \qquad \log A_n = X_n - \mu,$$

where $\mu > 0$ and $\{X_n\}$ satisfies (2.15). Then all of our hypotheses hold with $h = |f|$, provided that $\Lambda_B(\alpha) < \infty$ for some $\alpha > \mathfrak{r}$. The Markov chain $\{X_n\}$ is not uniformly recurrent, but the minorization $(\mathfrak{M})$ and other conditions are easily verified, as follows.

Let $P$ denote the transition kernel of $\{X_n\}$, and let $\Phi_x$ denote the Normal$(x, 1)$ density function. For any fixed $a > 0$, set

$$\underline{\Phi}_a(y) = \inf\{\Phi_{cx}(y) : |x| \le a\} = \min_{x \in \{-a,a\}} \Phi_{cx}(y) \qquad \forall y \in \mathbb{R}.$$

Then

$$(3.6) \qquad \int_E \underline{\Phi}_a(y)\,dy \le P(x, E) \qquad \forall x \in \{\tilde{x} \in \mathbb{S} : |\tilde{x}| \le a\},\ E \in \mathcal{S}.$$

Hence $(\mathcal{M})$ holds with $\nu_a(dy) = b\underline{\Phi}_a(y)\,dy$ and $b \in (0, \infty)$ a normalizing constant.

To verify the remaining conditions, first compute the cumulant generating function,

$$(3.7) \qquad \Lambda(\alpha) = \limsup_{n \to \infty} \frac{1}{n} \log \mathbb{E}[e^{\alpha S_n}] = -\alpha m + \frac{\sigma^2 \alpha^2}{2} \qquad \forall \alpha,$$

where $m$ and $\sigma$ are positive constants. The form of the function on the right-hand side is obtained by observing that $S_n := \log A_1 + \cdots + \log A_n$ is clearly normally distributed, and so the computation of the limit in (3.7) reduces to calculating its limiting normalized mean and variance. Now $\Lambda$ is finite everywhere and we have chosen $h = |f|$. Consequently the general remarks in Section 2 may be applied to obtain (H$_1$), (H$_2$) and (H$_3$).

If $\{X_n\}$ is an AR($p$) process with $p > 1$, then $\{X_n\}$ is not itself a Markov chain, but instead we may consider the Markov chain $\mathcal{X}_n = (X_{np}, \ldots, X_{n(p-1)+1})$. The Markov chain $\{X_n\}$ does *not* satisfy the minorization $(\mathfrak{M})$ [since we would need to take $k > 1$ in $(\mathcal{M}_0)$]. However, the chain $\{\mathcal{X}_n\}$ does satisfy the minorization $(\mathcal{M}_1)$ presented in Section 6.2 below, and hence the results of this paper may still be applied. For further discussion of this case and the verification of our conditions here, see Section 6.2 below. Finally, if $\{X_n\}$ is a general ARMA($p, q$) process, then a slight modification of our assumptions is necessary, but the basic approach still applies, in essence; see Section 6.2 for details.

Financial modeling by means of ARMA processes is now quite standard. In an insurance context, this forms the basis for the so-called Wilkie model; for an introduction in this setting, see Panjer (1998). A general introduction to ARMA models and their applications can be found, for example, in Box and Jenkins (1976) or Brockwell and Davis (1991).

EXAMPLE 3.3. Assume once again the Cramér–Lundberg model for the insurance business, but assume that the investments are split between a fixed proportion in a stock and bank, respectively. Suppose that the bank investment earns interest at a constant rate, say $r$, while the logarithmic stock returns, $\{R_n^*\}$, are modeled according to a stochastic volatility model, namely $R_n^* = \sigma_n \xi_n$, where $\{\xi_n\}$ is an i.i.d. Gaussian sequence of random variables, and $\{\log \sigma_n\}$ is rather arbitrary and may, for example, be taken to be an ARMA$(p,q)$ process, independent of $\{\xi_n\}$. Such a choice for $\{\log \sigma_n\}$ is fairly typical, and with this choice, the statistical properties of $\{R_n^*\}$ become mathematically tractable, leading to their popularity as an alternative to, say, GARCH models. For a futher description of these models and their statistical properties see, for example, Mikosch (2003) and Davis and Mikosch (2008a, 2008b). Recent developments can also be found, for example, in Shephard (2005).

Under the above assumptions,

$$(3.8) \qquad A_n = (p(1+r) + (1-p)e^{R_n^*})^{-1} \qquad \text{where } R_n^* = \sigma_n \xi_n,$$

for some constant $p \in (0,1)$, and thus $A_n$ is deterministically bounded from above, uniformly in $n$.

In this example, the Markov chain is $\{X_n\} = \{(\log \sigma_n, \xi_n)\} \subseteq \mathbb{R}^2$, and we may take $h(x) = \|x\|$ for all $x \in \mathbb{R}^2$. Note, in particular, that we would *not* want to choose $h = |f|$ since: (i) in order for (𝔐) to be satisfied, we would need $h$ to tend to infinity when $\|X_n\| \nearrow \infty$; and (ii) in order for (H$_2$) to be satisfied, we would need $\Lambda$ to be finite in a neighborhood of zero, but $\{f(X_n)\}$ has heavy tails along the negative axis. Thus, $h = |f|$ is not a suitable choice in this case.

Now, although the Markov chain $\{(\log \sigma_n, \xi_n)\}$ is two-dimensional, the analysis of it simplifies considerably due to the fact that $\{\xi_n\}$ is i.i.d. and $\{\log \sigma_n\}$ is independent of $\{\xi_n\}$. Indeed, since $\mathbb{E}[e^{\alpha \xi_1}] < \infty$ for all $\alpha$, and since the $A$-sequence is deterministically bounded from above, it is actually sufficient to verify (𝔐), (H$_2$) and (H$_3$) for the process $\{\log \sigma_n\}$ [using that $h(x) := \|x\| \le |x_1| + |x_2|$ in the verification of (H$_2$), and using (2.14) in place of (2.11) for this verification]. Moreover, since $\{\log \sigma_n\}$ is assumed to be an AR(1) process, these properties can be obtained just as in the previous example, while more general ARMA processes may also be considered, as discussed in Section 6.2 below.

The verification of (H$_1$) is, however, somewhat more complicated than in the previous example. Since the function $\Lambda$ is actually infinite along the negative axis, the reasoning leading to (2.10)—which gave the positivity of $\mathfrak{r}$ in (H$_1$)—is no longer applicable. To circumvent this difficulty, introduce the left-truncated sums

$$S_n^{(M)} = f^{(M)}(X_1) + \cdots + f^{(M)}(X_n),$$

$$\text{where } f^{(M)}(x) = f(x) \vee -M \text{ and } M \nearrow \infty.$$

Then $S_n^{(M)} \geq S_n$ for all $n$ and $M$, and therefore $\Lambda^{(M)}(\alpha) \geq \Lambda(\alpha)$ for all $\alpha > 0$, where $\Lambda^{(M)}$ is defined the same way as $\Lambda$, but with $\{S_n^{(M)}\}$ in place of $\{S_n\}$. Then $\mathbb{E}_\pi[f^{(M)}(X)] < 0$ for large $M$, and consequently the reasoning leading to (2.10) applies to $\{S_n^{(M)}\}$ and yields $(\Theta^{(M)})'(0) < 0$, for large $M$, where $\exp(-\Theta^{(M)})$ is the convergence parameter associated with $\{S_n^{(M)}\}$ [as defined in Nummelin (1984), page 27; see also Section 5 below]. In Proposition 5.1 and Remark 5.1 below, we will show that under $(\mathfrak{M})$, (H$_2$) and (H$_3$), the function $\Theta^{(M)}$ is convex and

$$\Theta^{(M)}(\mathfrak{r}_M) = \Lambda^{(M)}(\mathfrak{r}_M) = 0, \qquad \text{where } \mathfrak{r}_M = \sup\{\alpha : \Lambda^{(M)}(\alpha) \leq 0\}.$$

Therefore, $(\Theta^{(M)})'(0) < 0 \Longrightarrow \mathfrak{r}_M > 0$. Since $\Lambda^{(M)} \geq \Lambda$ for all $\alpha > 0$, we consequently obtain (H$_1$).

EXAMPLE 3.4.  Consider a GARCH$(1,1)$ process with Markov regime switching, namely,

$$(3.9) \quad R_n^* = \sigma_n \xi_n, \qquad \text{where } \sigma_n^2 = a_0 + b_1 \sigma_{n-1}^2 + a_1 (R_{n-1}^*)^2, \ n = 1, 2, \ldots,$$

and where $(a_0, a_1, b_1)$ is now viewed as a random vector which will typically oscillate between a finite number of states. As before, $\{\xi_n\}$ is an i.i.d. sequence of Gaussian random variables. [The following discussion applies equally well for an ARCH$(1)$ process in place of a GARCH$(1,1)$ process.] Such models have been studied rather extensively; for an introduction, see Lange and Rahbek (2008) and references therein. One motivation for considering this type of dependence is that it yields some of the statistical properties of long-range dependence. This viewpoint has been introduced in Mikosch and Stărică (2004). In any case, it is reasonable to assume that the parameters $(a_0, a_1, b_1)$ will change over time due, for example, to external economic factors, and this gives intuitive motivation for the dependence in $(a_0, a_1, b_1)$.

Under Markov regime switching, it is usually assumed that the observed parameter values at time $n$ are given by

$$\mathbf{a}_{X_n} := (a_0(X_n), a_1(X_n), b_1(X_n)) \in \mathbb{R}^3,$$

where $\{X_n\}$ is a $k$-state Markov chain which is exogenous and hence unobserved. Given the transition matrix of the finite-state chain, the extremes of this econometric process can then be analyzed by our methods, just as in Example 3.1.

An interesting variant would be to assume that $\{X_n\}$ is an ARMA process, driven by certain observed economic factors such as, for example, traded volume in the market, and then the extremes may be analyzed according to the methods in Example 3.2. It should be remarked in this example that if *both* $(a_1, b_1)$ and $a_0$ depend on $\{X_n\}$, then one needs a slight modification of the formulation of Section 2, as will be discussed below in Section 6.1.

To summarize, our conditions hold for a wide class of financial processes, such as regime-switching models, ARMA models and certain stochastic volatility models. Here we have developed the connection especially with ARMA models and with more complicated models which can be constructed from these, but it should be emphasized that our conditions are actually quite general. Indeed, our primary assumption is that the underlying Markov chain is light-tailed and therefore satisfies the usual moment conditions of modern large deviations theory, and for such processes, our results hold under very weak regularity conditions. It should, moreover, be pointed out that the investment process in Example 3.3 is actually *not* light-tailed—it is the driving Markovian process which is light-tailed, and for this reason our conditions can still be verified there.

Finally, we mention that the recurrence equations studied here are also broadly relevant for a wide class of related processes arising outside of the areas of financial and insurance mathematics. For some recent results along these lines see, for example, Gnedin (2007) or the survey article of Aldous and Bandyopadhyay (2005). Applications in the direction of computer science can be found, for example, in Neininger and Rüschendorf (2005) and references therein.

## 4. Proofs of the main theorems.

4.1. *Sketch of the proofs.* We begin with a brief sketch of the main ideas. Our starting point is the well-known regeneration lemma of Athreya and Ney (1978) and Nummelin (1978). This lemma asserts the existence of a sequence of random times, $T_0, T_1, T_2, \ldots$, such that the blocks

$$\{X_{T_{i-1}}, \ldots, X_{T_i - 1}\}, \qquad i = 0, 1, \ldots,$$

are independent for $i \geq 0$ and identically distributed for $i \geq 1$ (where $T_{-1} = 1$). Thus, in particular, the increments

$$(4.1) \qquad S_{T_i - 1} - S_{T_{i-1} - 1} := f(X_{T_{i-1}}) + \cdots + f(X_{T_i - 1})$$

are independent for $i \geq 0$ and identically distributed for $i \geq 1$. However, it is not immediately evident that a similar independence structure should exist for the risk process

$$(4.2) \qquad W_n := B_1 + A_1 B_2 + \cdots + (A_1 \cdots A_{n-1}) B_n,$$

which is our primary object of study. In Lemma 4.2 below, we shall show that

$$(4.3) \quad W_n = \breve{B}_0 + \breve{A}_0 \breve{B}_1 + \cdots + (\breve{A}_0 \cdots \breve{A}_{K^*-1}) \breve{B}_{K^*} + (\breve{A}_0 \cdots \breve{A}_{K^*}) \mathcal{R}_n,$$

where $\{(\breve{A}_i, \breve{B}_i) : i = 0, 1, \ldots\}$ is an independent sequence of random variables which, for $i \geq 1$, is identically distributed, $\mathcal{R}_n$ is a negligible remainder term, and $K^*(n) \nearrow \infty$ as $n \to \infty$. It follows as a consequence of (4.3) that $W := \sup_n W_n$ "nearly" satisfies the random recurrence equation

$$(4.4) \qquad W = \check{B} + \max\{\check{M} - \check{B}, \check{A}W\}$$

for some random variable $\check{M}$, where $(\check{A}, \check{B}) \overset{d}{=} (\check{A}_i, \check{B}_i)$ for $i \geq 1$. [A more precise statement will be given in Lemma 4.2.] The random recurrence equation (4.4) is very similar in form to (1.1), which applied in the i.i.d. case.

In Section 2 we argued that the probability of ruin is given by $\psi(u) := \mathbb{P}\{W > u\}$. Thus, to determine this probability, we need to find the tail distribution of $W$. To this end, we apply a result of Goldie (1991), which states that under appropriate conditions, (4.4) implies

$$(4.5) \qquad \lim_{u \to \infty} u^\eta \mathbb{P}\{W > u\} = D$$

for certain constants $D$ and $\eta$. Equation (4.5) provides a complete solution to our problem, but it remains to check that Goldie's conditions are actually satisfied and to identify the constant $\eta$ (which will be shown to equal $\mathfrak{r}$ in the proof of Theorem 2.1).

It is, in fact, quite challenging to verify that Goldie's conditions actually hold. To do so, we shall need to establish certain moment conditions for the random variables $\check{A}, \check{B}, \check{M}$ appearing in (4.4); in particular,

$$(4.6) \qquad \mathbb{E}[\check{A}^\alpha] < \infty, \qquad \mathbb{E}[|\check{B}|^\alpha] < \infty \quad \text{and} \quad \mathbb{E}[|\check{M}|^\alpha] < \infty$$

$$\text{for some } \alpha > \mathfrak{r};$$

see Theorem 4.2 below. The simplest of these studies $\check{A}^\alpha \overset{d}{=} \exp\{\alpha(S_{T_{i-1}} - S_{T_{i-1}})\}$; cf. (4.1). After a change of measure it can be shown that, roughly speaking,

$$(4.7) \qquad \mathbb{E}[\check{A}^\alpha] \approx \mathbb{E}^Q[e^{\varepsilon(T_i - T_{i-1})}] \qquad \text{for some } \varepsilon > 0,$$

for an $\alpha$-shifted kernel $Q$, and thus we see that such moment conditions are closely related to geometric recurrence for the $\alpha$-shifted chain. The random variables $|\check{B}|^\alpha$ and $|\check{M}|^\alpha$ are more complicated, but can be handled by somewhat similar techniques.

We now proceed more formally.

4.2. *Formal proofs of Theorems 2.1 and 2.2.* We begin with a precise statement of the regeneration lemma, first established by Athreya and Ney (1978) and Nummelin (1978). Here and in the following, set $T_{-1} = 1$ and let $\tau \overset{d}{=} T_i - T_{i-1}$ $(i \geq 1)$ denote a typical regeneration time.

Lemma 4.1. *Assume $(\mathcal{M}_0)$ holds with $k = 1$. Then there exists a sequence of random times, $0 < T_0 < T_1 < \cdots$, such that:*

(i)   $T_0, T_1 - T_0, T_2 - T_1, \ldots$ are finite a.s. and mutually independent;
(ii)  the sequence $\{T_i - T_{i-1} : i = 1, 2, \ldots\}$ is i.i.d.;
(iii) the random blocks $\{X_{T_{i-1}}, \ldots, X_{T_i-1}\}$ are independent, $i = 0, 1, \ldots$;
(iv)  $\mathbb{P}\{X_{T_i} \in E | \mathfrak{F}_{T_i-1}\} = \nu(E)$, for all $E \in \mathcal{S}$.

REMARK 4.1.   The regeneration times $\{T_i\}_{i \geq 0}$ of Lemma 4.1 can be related to the return times of $\{X_n\}$ to the set $C$ in $(\mathcal{M}_0)$, as follows. First introduce an augmented chain $\{(X_n, Y_n)\}$, where $\{Y_n : n = 0, 1, \ldots\}$ is an i.i.d. sequence of Bernoulli random variables, independent of $\{X_n\}$, with $\mathbb{P}\{Y_n = 1\} = \delta$, where $\delta \in (0, 1]$ is given as in $(\mathcal{M}_0)$. Then $T_i - 1$ can be identified as the $(i + 1)$th return time of $\{(X_n, Y_n)\}$ to the set $C \times \{1\}$. At the subsequent time, $T_i$, the random variable $X_{T_i}$ has the distribution $\nu$ given in $(\mathcal{M}_0)$, independent of the past history of the Markov chain.

To apply the lemma in the context of our problem, assume now that $(\mathcal{M})$ holds, and let $T_0, T_1, \ldots$ denote the resulting regeneration times. Then define the following random quantities formed over the independent random blocks of the previous lemma:

$$\check{A}_i = A_{T_{i-1}} \cdots A_{T_i-1}, \qquad i = 0, 1, \ldots;$$

$$\check{B}_i = B_{T_{i-1}} + A_{T_{i-1}} B_{T_{i-1}+1} + \cdots + (A_{T_{i-1}} \cdots A_{T_i-2}) B_{T_i-1}, \qquad i = 0, 1, \ldots;$$

$$\check{M}_i = \sup\{B_{T_{i-1}} + A_{T_{i-1}} B_{T_{i-1}+1} + \cdots + (A_{T_{i-1}} \cdots A_{j-1}) B_j : T_{i-1} \leq j < T_i\}, \qquad i = 0, 1, \ldots.$$

By Lemma 4.1, $\{(\check{A}_i, \check{B}_i, \check{M}_i) : i = 1, 2, \ldots\}$ is an i.i.d. sequence of random vectors which is also independent of $(\check{A}_0, \check{B}_0, \check{M}_0)$.

Our objective is to study the tail behavior of the random variable $W$ defined in Section 2.1. Let $W^{\mathfrak{R}}$ be defined in the same as $W$, but under the assumption that $T_0 = 1$, that is, under the assumption that regeneration occurs at time one. Also let $(\check{A}, \check{B}, \check{M}) \stackrel{d}{=} (\check{A}_1, \check{B}_1, \check{M}_1)$, and assume that $(\check{A}, \check{B}, \check{M})$ is independent of $\{(\check{A}_i, \check{B}_i, \check{M}_i) : i = 0, 1, \ldots\}$.

We begin by establishing the following.

LEMMA 4.2.   $W^{\mathfrak{R}} \stackrel{d}{=} \check{B} + \max\{\check{M} - \check{B}, \check{A} W^{\mathfrak{R}}\}$ and $W \stackrel{d}{=} \check{B}_0 + \max\{\check{M}_0 - \check{B}_0, \check{A}_0 W^{\mathfrak{R}}\}$.

PROOF.   For $n \geq 1$, set

$$K^* = K^*(n) := \inf\{i : T_i > n\} - 1.$$

Then for any $n \geq 1$,

$$W_n = (B_1 + A_1 B_2 + \cdots + (A_1 \cdots A_{T_0-2}) B_{T_0-1})$$

$$+ (A_1 \cdots A_{T_0 - 1})(B_{T_0} + \cdots + (A_{T_0} \cdots A_{T_1 - 2}) B_{T_1 - 1})$$

$$+ \cdots + (A_1 \cdots A_{T_{K^*-1}-1})(B_{T_{K^*-1}} + \cdots + (A_{T_{K^*-1}} \cdots A_{T_{K^*}-2}) B_{T_{K^*}-1})$$

$$+ (A_1 \cdots A_{T_{K^*}-1})(B_{T_{K^*}} + \cdots + (A_{T_{K^*}} \cdots A_{n-1}) B_n).$$

Hence

$$(4.8) \quad W_n = \check{B}_0 + \check{A}_0 \check{B}_1 + \cdots + (\check{A}_0 \cdots \check{A}_{K^*-1}) \check{B}_{K^*} + (\check{A}_0 \cdots \check{A}_{K^*}) \mathcal{R}_n,$$

where

$$\mathcal{R}_n := B_{T_{K^*}} + A_{T_{K^*}} B_{T_{K^*}+1} + \cdots + (A_{T_{K^*}} \cdots A_{n-1}) B_n.$$

Note that this last definition and the definition of $\check{M}_i$ imply

$$\check{M}_i = \sup\{\mathcal{R}_n : T_{i-1} \le n < T_i\} \qquad \text{for all } i = 0, 1, \ldots.$$

Hence by (4.8),

$$(4.9) \qquad\qquad W := \sup_{n \ge 1} W_n = \sup_{i \ge 0} V_i,$$

where

$$V_i := \begin{cases} \check{M}_0, & i = 0; \\ \check{B}_0 + \check{A}_0 \check{M}_1, & i = 1; \\ \check{B}_0 + \check{A}_0 \check{B}_1 + \cdots + (\check{A}_0 \cdots \check{A}_{i-2}) \check{B}_{i-1} \\ \quad + (\check{A}_0 \cdots \check{A}_{i-1}) \check{M}_i, & i = 2, 3, \ldots. \end{cases}$$

Moreover, by a repetition of the same argument,

$$(4.10) \qquad\qquad W^{\mathfrak{R}} = \sup_{i \ge 0} V_i^{(1)},$$

where $V_i^{(1)}$ is defined the same as $V_i$, but with all subscripts increased by a factor of one, so that $V_0^{(1)} = \check{M}_1$, and so on.

Next observe that

$$W = \sup_{i \ge 0} V_i := \sup\{\check{M}_0, \check{B}_0 + \check{A}_0 \check{M}_1, \check{B}_0 + \check{A}_0 \check{B}_1 + \check{A}_0 \check{A}_1 \check{M}_2, \ldots\}$$

$$(4.11) \qquad = \check{B}_0 + \sup\{\check{M}_0 - \check{B}_0, \check{A}_0 \check{M}_1, \check{A}_0(\check{B}_1 + \check{A}_1 \check{M}_2), \ldots\}$$

$$= \check{B}_0 + \max\{\check{M}_0 - \check{B}_0, \check{A}_0 W^{(1)}\},$$

where

$$W^{(j)} := \sup\{\check{M}_j, \check{B}_j + \check{A}_j \check{M}_{j+1}, \check{B}_j + \check{A}_j \check{M}_{j+1} + \check{A}_j \check{A}_{j+1} \check{M}_{j+2}, \ldots\} \overset{d}{=} W^{\mathfrak{R}}$$

$$\forall j \ge 1.$$

This establishes the second assertion of the lemma. For the first assertion, note by a repetition of (4.11) that

$$W^{\mathfrak{R}} = \sup_{i \geq 0} V_i^{(1)} = \check{B}_1 + \max\{\check{M}_1 - \check{B}_1, \check{A}_1 W^{(2)}\}$$

$$\overset{d}{=} \check{B} + \max\{\check{M} - \check{B}, \check{A} W^{\mathfrak{R}}\}. \qquad \square$$

Set

$$(4.12) \qquad \eta = \sup\{\alpha : \log \mathbb{E}[\check{A}^\alpha] \leq 0\}.$$

Now under our basic assumptions on $\{A_n\}$, $\log \check{A}$ is nonarithmetic. Hence by combining Lemma 4.2 with Theorem 6.2 of Goldie (1991), we obtain:

THEOREM 4.1. *Suppose $\eta > 0$ and $\mathbb{E}[\check{A}^\alpha] < \infty$ for some $\alpha > \eta$. Further assume*

$$(4.13) \qquad \mathbb{E}[|\check{B}|^\eta] < \infty \quad and \quad \mathbb{E}[|\check{M} - \check{B}|^\eta] < \infty.$$

*Then*

$$(4.14) \qquad \lim_{u \to \infty} u^\eta \mathbb{P}\{W^{\mathfrak{R}} > u\} = D,$$

*where $D \in [0, \infty)$ is given by*

$$(4.15) \qquad D = \frac{1}{\eta \tilde{m}} \mathbb{E}[((\check{B} + \max\{\check{M} - \check{B}, \check{A} W^{\mathfrak{R}}\})^+)^\eta - ((\check{A} W^{\mathfrak{R}})^+)^\eta]$$

*and $\tilde{m} := \mathbb{E}[\check{A}^\eta \log \check{A}]$.*

To apply the above result, we first need to develop some properties loosely related to geometric $\mathfrak{r}$-recurrence, which describe the moments associated with $\check{A}$, $\check{B}$ and $\check{M}$.

THEOREM 4.2. *Assume that $(\mathfrak{M})$, $(H_2)$ and $(H_3)$ are satisfied. Then for sufficiently large $a > 0$, there exists an $\alpha > \mathfrak{r}$ such that*

$$(4.16) \qquad \mathbb{E}[\check{A}^\alpha] < \infty, \qquad \mathbb{E}[|\check{B}|^\alpha] < \infty \quad and \quad \mathbb{E}[|\check{M}|^\alpha] < \infty$$

*and for a.a. $x \in \mathbb{S}$,*

$$(4.17) \qquad \mathbb{E}_x[\check{A}_0^\alpha] < \infty, \qquad \mathbb{E}_x[|\check{B}_0|^\alpha] < \infty \quad and \quad \mathbb{E}_x[|\check{M}_0|^\alpha] < \infty.$$

The proof of Theorem 4.2 is not straightforward and poses a main mathematical obstacle for our approach. This proof will be given below in Section 5. We now proceed directly to the proofs of Theorems 2.1 and 2.2.

PROOF OF THEOREM 2.1.    First assume $\mathfrak{r} = \eta$, where $\eta$ is defined as in (4.12). By Theorem 4.2, the conditions of Theorem 4.1 are then satisfied; and hence by (4.14) we obtain

$$(4.18) \qquad \lim_{u \to \infty} u^{\mathfrak{r}} \mathbb{P}\{W^{\mathfrak{R}} > u\} = D \qquad \text{as } u \to \infty,$$

where $D < \infty$ is given as in (4.15).

Since $\mathbb{E}[\check{A}^{\alpha}] < \infty$, some $\alpha > \mathfrak{r}$, it follows by (4.18) and a result of Breiman (1965) that

$$(4.19) \qquad \lim_{u \to \infty} u^{\mathfrak{r}} \mathbb{P}\{\check{A}_0 W^{\mathfrak{R}} > u\} = D \mathbb{E}[\check{A}_0^{\mathfrak{r}}] \qquad \text{as } u \to \infty.$$

Also, by Theorem 4.2 and an application of Chebyshev's inequality,

$$(4.20) \qquad \mathbb{P}\{\check{B}_0 > u\} \le D_1 u^{-\alpha} \quad \text{and} \quad \mathbb{P}\{\check{M}_0 - \check{B}_0 > u\} \le D_2 u^{-\alpha},$$

where $\alpha > \mathfrak{r}$ and $D_1$, $D_2$ are finite constants. Moreover by Lemma 4.2,

$$(4.21) \qquad W \overset{d}{=} \check{B}_0 + \max\{\check{M}_0 - \check{B}_0, \check{A}_0 W^{\mathfrak{R}}\}.$$

Hence

$$(4.22) \qquad \lim_{u \to \infty} u^{\mathfrak{r}} \mathbb{P}\{W > u\} = D \mathbb{E}[\check{A}_0^{\mathfrak{r}}] \qquad \text{as } u \to \infty,$$

which establishes (2.20).

It remains to show that $\mathfrak{r} = \eta$. Suppose false. Let $(\theta(\alpha))^{-1}$ denote the convergence parameter of the kernel $\hat{P}_\alpha$ in (2.8), and let $\Theta = \log \theta$. [For the definition, see Nummelin (1984), page 27.] Then it is well known that $\Theta(\alpha) \le \Lambda(\alpha)$; cf. Proposition 5.1(ii) below. Hence

$$(4.23) \qquad \mathfrak{r} = \sup\{\alpha : \Lambda(\alpha) \le 0\} \le \sup\{\alpha : \Theta(\alpha) \le 0\}.$$

Moreover, it follows from Nummelin [(1984), Proposition 4.7(ii)] that

$$(4.24) \qquad \mathbb{E}_{\nu_a}[e^{\alpha \check{S} - \tau \Theta(\alpha)}] \le 1,$$

where $\tau \overset{d}{=} T_i - T_{i-1}$ and $\check{S} \overset{d}{=} S_{T_i - 1} - S_{T_{i-1} - 1}$ for $i > 1$. [See Ney and Nummelin (1987a), Section 4, and Proposition 5.1(iv) below for closely related results.] It follows as a consequence of (4.24) that

$$(4.25) \qquad \Theta(\alpha) \le 0 \quad \Longrightarrow \quad \Pi(\alpha) := \mathbb{E}[e^{\alpha \check{S}}] \le 1$$

and therefore by (4.23), $\mathfrak{r} \le \eta$. Hence, if $\mathfrak{r} \ne \eta$, then we must have $\mathfrak{r} < \alpha < \eta$ for some $\alpha$.

Now assume that this is the case, and choose $\alpha \in (\mathfrak{r}, \eta)$. Set

$$L_i = \sup\{S_j - S_{T_{i-1}} : T_{i-1} \le j < T_i\}, \qquad i = 0, 1, \dots$$

(where $T_{-1} := 1$), and set $L \overset{d}{=} L_i$ for $i \geq 1$. By taking $\{B_1, B_2, \ldots\} = \{1, 1, \ldots\}$ in Theorem 4.2, we obtain that for sufficiently small $\alpha > \mathfrak{r}$,

$$(4.26) \quad \mathbb{E}[e^{\alpha L}] < \infty, \qquad \mathbb{E}[e^{\alpha L_0}] < \infty \quad \text{and} \quad \mathbb{E}[e^{\alpha S_{T_0 - 1}}] := \mathbb{E}[\check{A}_0^{\alpha}] < \infty.$$

Then by the independence of the regeneration cycles [cf. Lemma 4.1],

$$
\begin{aligned}
(4.27) \quad \mathbb{E}[e^{\alpha S_n}] &\leq \mathbb{E}\left[e^{\alpha L_0} + \sum_{i=1}^{\infty} \exp\{\alpha(S_{T_{i-1}-1} + L_i)\}\right] \\
&\leq \mathbb{E}[e^{\alpha L_0}] + \mathbb{E}[e^{\alpha L}]\mathbb{E}[e^{\alpha S_{T_0-1}}]\sum_{i=0}^{\infty}\mathbb{E}[e^{\alpha \check{S}}]^i < \infty,
\end{aligned}
$$

where the last step follows since $\alpha < \eta \implies \Pi(\alpha) < 1$, by the strict convexity of $\Pi$. By (4.27) and the definition of $\Lambda$, we conclude that $\Lambda(\alpha) \leq 0$, that is, $\alpha \leq \mathfrak{r}$, a contradiction. Therefore $\mathfrak{r} = \eta$. $\quad\square$

PROOF OF THEOREM 2.2. The proof is very similar to that of Theorem 2.1, but easier, so we only sketch the details. The main modification is in Lemma 4.2. In particular, we need to develop analogous random recurrence equations for $W_{\infty}$ and $W_{\infty}^*$.

Let $\{(\check{A}_i, \check{B}_i, \check{M}_i) : i = 0, 1, \ldots\}$ be defined as in the discussion prior to Lemma 4.2, and recall from the proof of Lemma 4.2 [cf. (4.8)] that

$$(4.28) \quad W_n = \check{B}_0 + \check{A}_0 \check{B}_1 + \cdots + (\check{A}_0 \cdots \check{A}_{K^*-1})\check{B}_{K^*} + (\check{A}_0 \cdots \check{A}_{K^*})\mathcal{R}_n,$$

where $K^* := \inf\{i : T_i > n\} - 1$.

Observe that $W_{\infty} := \lim_{n \to \infty} W_n$ exists a.s. by Theorem 2.1 of Goldie and Maller (2000) and Theorem 4.2 above. Moreover, since $K^* \to \infty$ a.s. as $n \to \infty$, the last term in (4.28) converges to zero a.s. as $n \to \infty$. Hence

$$(4.29) \quad W_{\infty} \overset{d}{=} \check{B}_0 + \check{A}_0 W_{\infty}^{\mathfrak{R}},$$

where $W_{\infty}^{\mathfrak{R}}$ is defined the same as $W_{\infty}$, except that now $T_0 = 1$, so that regeneration occurs at the initial time. Furthermore, a repetition of the argument leading to (4.29) yields

$$(4.30) \quad W_{\infty}^{\mathfrak{R}} \overset{d}{=} \check{B} + \check{A} W_{\infty}^{\mathfrak{R}},$$

where $(\check{A}, \check{B}) \overset{d}{=} (\check{A}_i, \check{B}_i), i \geq 1$; cf. the proof of Lemma 4.2. The last equation is a random recurrence equation, which can be analyzed using Goldie's (1991) Theorem 4.1. First note by Theorem 4.2 that the moment conditions in Goldie's theorem are satisfied. The required result then follows from (4.29), (4.30) and Goldie (1991), Theorem 4.1. [Specifically, by employing (4.29), (4.30) in place of Lemma 4.2 above, and by employing Goldie's Theorem 4.1 in place of his Theorem 6.2, and then reasoning as in the proof of Theorem 2.1, the desired limit result is obtained.]

For the process $\{W_n^*\}$, define $\{\check{A}_i\}$ as before, and let

$$\check{B}_i^* = (A_{T_{i-1}} \cdots A_{T_{i-1}+1})B_{T_{i-1}} + \cdots + A_{T_{i-1}}B_{T_{i-2}} + B_{T_{i-1}} \qquad \forall i \geq 0,$$

$$\check{M}_i^* = \sup\{(A_{T_{i-1}} \cdots A_{j+1})|B_j| + \cdots + |B_{T_{i-1}}| : T_{i-1} \leq j < T_i\} \qquad \forall i \geq 0.$$

[Note that the definition of $\check{B}_i^*$ differs from that of $\check{B}_i$, because a given element $B_n$ is now multiplied by $(A_{n+1} \cdots A_{T_{i-1}})$ rather than by $(A_{T_{i-1}} \cdots A_{n-1})$.]

Let $Z_i^* := W_{T_i-1}^*$ denote the value of $\{W_n^*\}$ at its $(i+1)$th regeneration time $(i = 0, 1, \dots)$. Since

$$W_n^* = B_n + A_n W_{n-1}^*, \qquad n = 1, 2, \dots,$$

it follows after a short argument that

$$(4.31) \qquad Z_i^* = \check{B}_i^* + \check{A}_i Z_{i-1}^*, \qquad i = 1, 2, \dots.$$

First consider $Z_\infty^* := \lim_{i \to \infty} Z_i^*$. By Theorem 3.1 of Goldie and Maller (2000) and Theorem 4.2 above, the limit exists a.s., and by (4.31), it satisfies the recurrence equation

$$(4.32) \qquad Z_\infty^* \overset{d}{=} \check{B}^* + \check{A} Z_\infty^*,$$

where $(\check{A}, \check{B}^*) \overset{d}{=} (\check{A}_i, \check{B}_i^*)$ for $i \geq 1$. Then by (4.32) and Goldie [(1991), Theorem 4.1],

$$(4.33) \qquad \mathbb{P}\{Z_\infty^* > u\} \sim \tilde{C} u^{-\mathfrak{r}} \qquad \text{as } u \to \infty$$

for some constant $\tilde{C} < \infty$. As before, the required moment conditions are verified using Theorem 4.2.

It remains to show that $\lim_{n \to \infty} W_n^* = \lim_{i \to \infty} Z_i^*$. To this end, let $L(n) := \inf\{i : T_i > n\}$ denote the first regeneration time after time $n$. Then, for example, the term "$B_n$" forms a part of the sum $\check{B}_{L(n)}^*$, and the discrepency of $W_n^*$ from $Z_{L(n)}^*$ is bounded by $\check{M}_{L(n)}$. More precisely,

$$(4.34) \qquad Z_{L(n)}^* - \check{M}_{L(n)}^* \leq W_n^* \leq Z_{L(n)}^* + \check{M}_{L(n)}^*.$$

By a minor variant of Theorem 4.2, $\check{M}^*$ is asymptotically negligible compared with $\lim_{n \to \infty} Z_{L(n)}^*$. Consequently $W_\infty^* = Z_\infty^*$. $\square$

## 5. Proof of Theorem 4.2 and some related regularity results.

5.1. *Notation and preliminary remarks.* In the following discussion, it will always be assumed that the minorization ($\mathcal{M}$) holds for a given parameter $a > 0$.

First we introduce some additional notation. For any kernel $K$, let

$$\hat{K}_\alpha(x, E) = \int_E e^{\alpha f(y)} K(x, dy) \qquad \forall \alpha \in \mathbb{R},$$

and let

$$\hat{K}_{\boldsymbol{\alpha}}(x,E) = \int_E e^{\alpha f(y)+\beta h(y)} K(x,dy) \qquad \forall \boldsymbol{\alpha} = (\alpha,\beta) \in \mathbb{R}^2.$$

Let $(\theta(\alpha))^{-1}$ denote the *convergence parameter* of $\hat{P}_\alpha$ [as defined in Nummelin (1984), page 27], and let $(\boldsymbol{\theta}(\boldsymbol{\alpha}))^{-1}$ denote the convergence parameter of $\hat{P}_{\boldsymbol{\alpha}}$. Let $\Theta(\alpha) = \log\theta(\alpha)$ and $\boldsymbol{\Theta}(\boldsymbol{\alpha}) = \log\boldsymbol{\theta}(\boldsymbol{\alpha})$, for all $\boldsymbol{\alpha} = (\alpha,\beta) \in \mathbb{R}^2$. [We will repeatedly use the notation $\boldsymbol{\alpha} = (\alpha,\beta)$ and generally do this without stating so explicitly.]

Below, we will often need to work with a perturbed kernel, namely,

$$P^{(\varepsilon)}(x,E) := P(x,E) + \varepsilon\nu_a(E) \qquad \forall E \in \mathcal{S}, \ \varepsilon \geq 0,$$

where $\nu_a$ is given as in $(\mathcal{M})$. In this connection, let $\hat{P}_\alpha^{(\varepsilon)}$, $\hat{P}_{\boldsymbol{\alpha}}^{(\varepsilon)}$, $\theta^{(\varepsilon)}$, $\boldsymbol{\theta}^{(\varepsilon)}$, $\Theta^{(\varepsilon)}$ and $\boldsymbol{\Theta}^{(\varepsilon)}$ be defined the same as $\hat{P}_\alpha$, $\hat{P}_{\boldsymbol{\alpha}}$ and so on, but with $P$ replaced everywhere in these definitions with $P^{(\varepsilon)}$.

For each $x \in \mathbb{S}$, set

$$r_\alpha^{(\varepsilon)}(x) = \sum_{n=0}^\infty (\theta^{(\varepsilon)}(\alpha))^{-n-1} (\hat{P}_\alpha^{(\varepsilon)} - \delta_a \mathbf{1}_{\mathcal{L}_a h} \otimes \hat{\nu}_a^{(\alpha)})^n \delta_a \mathbf{1}_{\mathcal{L}_a h}(x)$$

and

$$\mathbf{r}_{\boldsymbol{\alpha}}^{(\varepsilon)}(x) = \sum_{n=0}^\infty (\boldsymbol{\theta}^{(\varepsilon)}(\boldsymbol{\alpha}))^{-n-1} (\hat{P}_{\boldsymbol{\alpha}}^{(\varepsilon)} - \delta_a \mathbf{1}_{\mathcal{L}_a h} \otimes \hat{\nu}_a^{(\boldsymbol{\alpha})})^n \delta_a \mathbf{1}_{\mathcal{L}_a h}(x),$$

where, for arbitrary $g : \mathbb{S} \to \mathbb{R}$ and $\mu : \mathcal{S} \to \mathbb{R}$:

$$g \otimes \mu(x,E) = g(x)\mu(E); \qquad \hat{\mu}^{(\alpha)}(E) = \int_E e^{\alpha f(y)} \mu(dy);$$

$$\hat{\mu}^{(\boldsymbol{\alpha})}(E) = \int_E e^{\alpha f(y)+\beta h(y)} \mu(dy).$$

In the special case that $\varepsilon = 0$, we shall simply write $r_\alpha$, $\mathbf{r}_{\boldsymbol{\alpha}}$ in place of $r_\alpha^{(0)}$, $\mathbf{r}_{\boldsymbol{\alpha}}^{(0)}$.

The function $r_\alpha^{(\varepsilon)}$ is known to be $(\theta^{(\varepsilon)}(\alpha))^{-1}$-subinvariant with respect to the kernel $\hat{P}_\alpha$, and moreover to be $(\theta^{(\varepsilon)}(\alpha))^{-1}$-invariant in the case that $\hat{P}_\alpha$ is $(\theta^{(\varepsilon)}(\alpha))^{-1}$-recurrent; see Nummelin (1984), Theorem 5.1 and its proof [and the proof of Proposition 5.1(iv) below]. Likewise, the function $\mathbf{r}_{\boldsymbol{\alpha}}^{(\varepsilon)}$ is known to be $(\boldsymbol{\theta}^{(\varepsilon)}(\boldsymbol{\alpha}))^{-1}$-subinvariant with respect to the kernel $\hat{P}_{\boldsymbol{\alpha}}$, and to be $(\boldsymbol{\theta}^{(\varepsilon)}(\boldsymbol{\alpha}))^{-1}$-invariant in the case that $\hat{P}_{\boldsymbol{\alpha}}$ is $(\boldsymbol{\theta}^{(\varepsilon)}(\boldsymbol{\alpha}))^{-1}$-recurrent.

Let

$$\boldsymbol{\Lambda}^{(\varepsilon)}(\boldsymbol{\alpha}) = \limsup_{n\to\infty} \frac{1}{n} \log(\hat{P}_{\boldsymbol{\alpha}}^{(\varepsilon)})^n(x,\mathbb{S}) \qquad \forall \boldsymbol{\alpha} \in \mathbb{R}^2,$$

where $x$ is the initial state of the Markov chain. [If we replace $P^{(\varepsilon)}$ with $P$ in this definition, then the right-hand side reduces to the Gärtner–Ellis limit, $\boldsymbol{\Lambda}$, which was introduced in Section 2.]

Finally, in the proofs below we will make use of the "shifted" kernel

$$Q_\alpha^{(\varepsilon)}(x, E) := \int_E \frac{e^{\alpha f(y)} r_\alpha^{(\varepsilon)}(y)}{\theta^{(\varepsilon)}(\alpha) r_\alpha^{(\varepsilon)}(x)} P^{(\varepsilon)}(x, dy).$$

Under the minorization $(\mathcal{M})$, observe that the kernel $Q_\alpha^{(\varepsilon)}$ itself satisfies a minorization, namely,

$(\mathcal{M}_Q)$ $\qquad g_\alpha^{(\varepsilon)}(x) \mu_\alpha^{(\varepsilon)}(E) \le Q_\alpha^{(\varepsilon)}(x, E) \qquad \forall x \in \mathbb{S}, \; E \in \mathcal{S},$

where

$$g_\alpha^{(\varepsilon)}(x) := \left( \frac{\delta_a B}{\theta^{(\varepsilon)}(\alpha) r_\alpha^{(\varepsilon)}(x)} \mathbf{1}_{\mathcal{L}_a h}(x) \right) \wedge \frac{1}{2},$$

$$\mu_\alpha^{(\varepsilon)}(E) := \frac{1}{B} \int_E e^{\alpha f(y)} r_\alpha^{(\varepsilon)}(y) \nu_a(dy)$$

and $B$ is a normalizing constant chosen so that $\mu_\alpha^{(\varepsilon)}$ is a probability measure. After a truncation of $\nu_a$ as described above in Remark 2.1, the integral in the definition of $\mu_\alpha^{(\varepsilon)}$ will always be finite and hence $0 < B < \infty$.

5.2. *Some regularity properties.* In the next two propositions, we collect various regularity properties which will be needed in Section 5.3. After reading the statement of the propositions, the reader may want to proceed directly to Section 5.3—which contains the core results of this part of the paper—and refer back to the present section as necessary.

In the following discussion, set $\inf\{\varnothing\} = \infty$ and let $\mathfrak{D}_\Lambda$, $\mathfrak{D}_{\boldsymbol{\Lambda}}$ and $\mathfrak{D}_{\boldsymbol{\Theta}}$ denote the domains of $\Lambda$, $\boldsymbol{\Lambda}$ and $\boldsymbol{\Theta}$, respectively. [To be entirely precise, $\boldsymbol{\alpha} \in \mathfrak{D}_{\boldsymbol{\Lambda}}$ means that $\boldsymbol{\Lambda}(\boldsymbol{\alpha}) < \infty$ for $\varphi$-a.a. initial states $x$, cf. part (i) of the next proposition; and our main results will be valid away from an appropriate set of measure zero.]

We remark that in (i)–(iii) of the following proposition, we develop properties of $\boldsymbol{\Lambda}^{(\varepsilon)}(\boldsymbol{\alpha})$, $\boldsymbol{\Theta}^{(\varepsilon)}(\boldsymbol{\alpha})$, etc., but the same properties hold also for $\Lambda^{(\varepsilon)}(\alpha)$, $\Theta^{(\varepsilon)}(\alpha)$, etc., as can be seen by setting $\boldsymbol{\alpha} = (\alpha, 0)$ and observing that $\boldsymbol{\Lambda}^{(\varepsilon)}(\boldsymbol{\alpha}) = \Lambda^{(\varepsilon)}(\alpha)$, $\boldsymbol{\Theta}^{(\varepsilon)}(\boldsymbol{\alpha}) = \Theta^{(\varepsilon)}(\alpha)$, and so on.

PROPOSITION 5.1. *Assume* $(\mathcal{M})$. *Then:*

(i) *For any* $\boldsymbol{\alpha} \in \mathfrak{D}_{\boldsymbol{\Lambda}}$,

(5.1) $\qquad\qquad \boldsymbol{\Lambda}(\boldsymbol{\alpha}) = \boldsymbol{\Lambda}_{\nu_a}(\boldsymbol{\alpha}), \qquad \varphi\text{-}a.a.\ x,$

*where $\boldsymbol{\Lambda}_{\nu_a}$ denotes the same limiting quantity as $\boldsymbol{\Lambda}$, but conditioned on regeneration at time zero. Thus, $\boldsymbol{\Lambda}$ is independent of its initial state. Moreover for all $N \geq 1$,*

$$(5.2) \qquad \bar{\boldsymbol{\Lambda}}_N(\boldsymbol{\alpha}) := \sup_{n \geq N} \frac{1}{n} \log \mathbb{E}_{\nu_a}^{(1)}[e^{\alpha S_n + \beta S_n^{(h)}}] < \infty \qquad \forall \boldsymbol{\alpha} \in \mathfrak{D}_{\boldsymbol{\Lambda}},$$

*where $\mathbb{E}_{\nu_a}^{(1)}[\cdot]$ denotes that regeneration occurs at time one.*

(ii) *The function $\boldsymbol{\Theta}^{(\varepsilon)}$ is convex, and $\boldsymbol{\Theta}^{(\varepsilon)}(\boldsymbol{\alpha}) \leq \boldsymbol{\Lambda}^{(\varepsilon)}(\boldsymbol{\alpha})$ for all $\boldsymbol{\alpha} \in \mathfrak{D}_{\boldsymbol{\Lambda}}$ and $\varepsilon \geq 0$.*

(iii) *For any $\boldsymbol{\alpha} \in \mathfrak{D}_{\boldsymbol{\Lambda}}$,*

$$(5.3) \qquad \lim_{\varepsilon \to 0} \boldsymbol{\Lambda}^{(\varepsilon)}(\boldsymbol{\alpha}) = \boldsymbol{\Lambda}(\boldsymbol{\alpha}).$$

(iv) *If $\alpha \in \mathfrak{D}_{\Lambda}$ and $\varepsilon \geq 0$, then $Q_\alpha^{(\varepsilon)}$ is either a subprobability or probability measure. Moreover, if $\tilde{P}_\alpha^{(\varepsilon)}$ is $(\theta(\alpha))^{-1}$-recurrent, then $Q_\alpha^{(\varepsilon)}$ is actually a probability measure.*

PROOF. (i) Let $\Delta > 0$, and set

$$F_\Delta = \{x \in \mathbb{S} \colon \boldsymbol{\Lambda}_x(\boldsymbol{\alpha}) \geq \boldsymbol{\Lambda}_{\nu_a}(\boldsymbol{\alpha}) + \Delta\},$$

where $\boldsymbol{\Lambda}_x$ denotes—now explicitly—that we are conditioning on the Markov chain starting in state $x$. Also, set

$$\mu_m(E) = \mathbb{E}_{\nu_a}[e^{\langle \boldsymbol{\alpha}, \mathbf{S}_m \rangle}; X_m \in E] \qquad \forall m \in \mathbb{Z}_+, \ E \in \mathcal{S},$$

where $\mathbf{S}_n = (S_n, S_n^{(h)})$.

Note that

$$(5.4) \qquad \mathbb{E}_{\nu_a}[e^{\langle \boldsymbol{\alpha}, \mathbf{S}_{m+n} \rangle}] \geq \mathbb{E}_{\nu_a}[e^{\langle \boldsymbol{\alpha}, \mathbf{S}_m \rangle}; X_m \in F_\Delta] \inf_{x \in F_\Delta} \mathbb{E}_x[e^{\langle \boldsymbol{\alpha}, \mathbf{S}_n \rangle}].$$

Hence it follows from the definitions of $F_\Delta$ and $\{\mu_m\}$ that

$$(5.5) \qquad \boldsymbol{\Lambda}_{\nu_a}(\boldsymbol{\alpha}) \geq \liminf_{n \to \infty} \frac{1}{n} \log \mu_m(F_\Delta) + (\boldsymbol{\Lambda}_{\nu_a}(\boldsymbol{\alpha}) + \Delta).$$

Since $\mu_m(F_\Delta)$ obviously does not depend on $n$, this last equation is only possible if $\mu_m(F_\Delta) = 0$ for all $m$. Consequently, $\mathbb{P}\{X_m \in F_\Delta\} = 0$ for all $m$. Since $\varphi(F_\Delta) > 0$ would imply that $\{X_m\}$ would visit $F_\Delta$ with positive probability over a regeneration cycle [Athreya and Ney (1978), Section 6], we conclude that $\varphi(F_\Delta) = 0$. Since the last equality holds for any $\Delta > 0$, it now follows from the definition of $F_\Delta$ that $\boldsymbol{\Lambda}_x(\boldsymbol{\alpha}) \leq \boldsymbol{\Lambda}_{\nu_a}(\boldsymbol{\alpha})$ for $\varphi$-a.a. $x$.

Conversely, if $\mathcal{E}_m$ denotes the event that regeneration occurs at time $m$, then

$$(5.6) \qquad \mathbb{E}_x[e^{\langle \boldsymbol{\alpha}, \mathbf{S}_{m+n} \rangle}] \geq \mathbb{E}_x[e^{\langle \boldsymbol{\alpha}, \mathbf{S}_m \rangle}; \mathcal{E}_m] \mathbb{E}_{\nu_a}[e^{\langle \boldsymbol{\alpha}, \mathbf{S}_n \rangle}].$$

If $m$ is chosen such that $\mathbb{P}\{\mathcal{E}_m\} > 0$, then $\mathbb{E}_x[e^{\langle \boldsymbol{\alpha}, \mathbf{S}_m \rangle}; \mathcal{E}_m] > 0$. Hence it follows upon taking $\limsup_{n \to \infty} n^{-1} \log(\cdot)$ in (5.6) that $\boldsymbol{\Lambda}_x(\boldsymbol{\alpha}) \geq \boldsymbol{\Lambda}_{\nu_a}(\boldsymbol{\alpha})$.

It remains to study (5.2). Taking $\limsup_{m \to \infty} m^{-1} \log(\cdot)$ on the left- and right-hand sides of (5.6), where $\mathcal{E}_m$ now denotes the event that $X_m \in \mathcal{L}_a h$ [the "small set" in $(\mathcal{M})$], yields

$$(5.7) \qquad \boldsymbol{\Lambda}(\boldsymbol{\alpha}) \geq \boldsymbol{\Theta}(\boldsymbol{\alpha}) + \lim_{m \to \infty} \frac{1}{m} \log \mathbb{E}_{\nu_a}^{(1)}[e^{\langle \boldsymbol{\alpha}, \mathbf{S}_n \rangle}],$$

where the first term on the right-hand side was obtained from the definition of the convergence parameter [cf. Nummelin (1984), page 27. On the right-hand side, we have also used the fact that regeneration occurs with probability $\delta_a$ upon each return to $\mathcal{L}_a h$, independent of the prior evolution of the chain].

Now the term on the left-hand side of (5.7) is finite, by assumption; while the first term on the right-hand side is greater than $-\infty$, by Nummelin (1984), Theorem 3.2. Since the second term on the right (but inside the limit) does not depend on $m$, we conclude that it must also be finite, that is,

$$(5.8) \qquad \mathbb{E}_{\nu_a}^{(1)}[e^{\langle \boldsymbol{\alpha}, \mathbf{S}_n \rangle}] < \infty \qquad \forall n \in \mathbb{Z}_+, \ \boldsymbol{\alpha} \in \mathfrak{D}_{\boldsymbol{\Lambda}}.$$

Furthermore, for any $\boldsymbol{\alpha} \in \mathfrak{D}_{\boldsymbol{\Lambda}}$, it follows by Remark 2.1 (and the fact that we have just shown $\boldsymbol{\Lambda} = \boldsymbol{\Lambda}_{\nu_a}$) that

$$(5.9) \qquad \limsup_{n \to \infty} \frac{1}{n} \log \mathbb{E}_{\nu_a}^{(1)}[e^{\langle \boldsymbol{\alpha}, \mathbf{S}_n \rangle}] = \boldsymbol{\Lambda}(\boldsymbol{\alpha}) < \infty.$$

[After truncation of $\nu_a$, as described in Remark 2.1, we have $\mathrm{supp}\,\nu_a \subseteq \mathcal{L}_a h \cap \mathcal{L}_b f$, where $\max\{a, b\} < \infty$, and hence $e^{\langle \boldsymbol{\alpha}, \mathbf{S}_1 \rangle}$ is deterministically bounded, meaning that the superscript "(1)" may be removed when taking the limit on the left-hand side of (5.9).] Hence, the terms inside the limit in (5.9) are all finite and, in the limit, they converge to a finite constant. Consequently, (5.2) follows from (5.8) and (5.9).

(ii) The convexity of $\boldsymbol{\Theta}^{(\varepsilon)}$ is obtained as in Iscoe, Ney and Nummelin (1985), Lemma 3.4. For further details, see Collamore (2002), Lemma 4.2.

To see that $\boldsymbol{\Theta}^{(\varepsilon)} \leq \boldsymbol{\Lambda}^{(\varepsilon)}$, note that the minorization $(\mathcal{M})$ implies a corresponding minorization for $\hat{P}_{\boldsymbol{\alpha}}$, namely,

$$(\hat{\mathcal{M}}) \qquad \delta_a \mathbf{1}_{\mathcal{L}_a h}(x) \hat{\nu}_a^{(\boldsymbol{\alpha})}(E) \leq \hat{P}_{\boldsymbol{\alpha}}(x, E) \qquad \forall x \in \mathbb{S}, \ E \in \mathcal{S}.$$

Then by the definition of the convergence parameter [Nummelin (1984), page 27] and Proposition 3.4 of Nummelin (1984), we see that $(\boldsymbol{\Theta}^{(\varepsilon)}(\boldsymbol{\alpha}))^{-1}$ is the radius of convergence of the power series

$$\sum_{n=0}^{\infty} \gamma^n \hat{\nu}_a^{(\boldsymbol{\alpha})}(\hat{P}_{\boldsymbol{\alpha}}^{(\varepsilon)})^n(\mathcal{L}_a h)$$

(since $\mathcal{L}_a h$ is a "small set"). But $(\mathbf{\Lambda}^{(\varepsilon)}(\boldsymbol{\alpha}))^{-1}$ is the radius of convergence of

$$\sum_{n=0}^{\infty} \gamma^n \hat{\nu}_a^{(\boldsymbol{\alpha})} (\hat{P}_{\boldsymbol{\alpha}}^{(\varepsilon)})^n (\mathbb{S}).$$

Hence $\mathcal{L}_a h \subseteq \mathbb{S} \Longrightarrow \mathbf{\Theta}^{(\varepsilon)}(\boldsymbol{\alpha}) \le \mathbf{\Lambda}^{(\varepsilon)}(\boldsymbol{\alpha})$.

(iii) First observe [following the proof in Collamore (2002), Theorem 3.1] that

$$
\begin{aligned}
\mathfrak{L}_n^{(\varepsilon)}(\boldsymbol{\alpha}) &:= \int_{\mathbb{S}} \hat{\nu}_a^{(\boldsymbol{\alpha})}(dx) (\hat{P}_{\boldsymbol{\alpha}}^{(\varepsilon)})^{n-1}(x, \mathbb{S}) \\
(5.10) \qquad &= \sum_{(i_1, \ldots, i_{n-1}) \in \mathbb{I}} \int_{\mathbb{S}^n} \hat{\nu}_a^{(\boldsymbol{\alpha})}(dx_1) J^{(i_1)}(x_1, dx_2) \cdots J^{(i_{n-1})}(x_{n-1}, dx_n),
\end{aligned}
$$

where $J^{(0)} = \varepsilon \hat{\nu}_a$, $J^{(1)} = \hat{P}_{\boldsymbol{\alpha}}$ and $\mathbb{I}$ consists of all elements of the form $(i_1, \ldots, i_{n-1})$ such that $i_j \in \{0, 1\}$, $j = 1, \ldots, n-1$. What needs to be shown is that

$$(5.11) \qquad \mathbf{\Lambda}_{\nu_a}^{(\varepsilon)}(\boldsymbol{\alpha}) := \limsup_{n \to \infty} \frac{1}{n} \log \mathfrak{L}_n^{(\varepsilon)}(\boldsymbol{\alpha}) \to \mathbf{\Lambda}_{\nu_a}(\boldsymbol{\alpha}) \qquad \text{as } \varepsilon \to 0.$$

[The required result then follows from (5.11), since by (i) we have $\mathbf{\Lambda}(\boldsymbol{\alpha}) = \mathbf{\Lambda}_{\nu_a}(\boldsymbol{\alpha})$, and a repetition of the same argument also gives $\mathbf{\Lambda}^{(\varepsilon)}(\boldsymbol{\alpha}) = \mathbf{\Lambda}_{\nu_a}^{(\varepsilon)}(\boldsymbol{\alpha})$.]

Let

$$\log \mathfrak{b}_N = N\{\bar{\mathbf{\Lambda}}_1(\boldsymbol{\alpha}) - \bar{\mathbf{\Lambda}}_N(\boldsymbol{\alpha})\} \qquad \forall N \ge 1,$$

where $\bar{\mathbf{\Lambda}}_N$ is defined as in (5.2). Setting $\sigma_N(\boldsymbol{\alpha}) = \exp \bar{\mathbf{\Lambda}}_N(\boldsymbol{\alpha})$, then by definition

$$(5.12) \qquad \mathfrak{b}_N = \left(\frac{\sigma_1(\boldsymbol{\alpha})}{\sigma_N(\boldsymbol{\alpha})}\right)^N.$$

Note that if $\boldsymbol{\alpha} \in \mathfrak{D}_{\mathbf{\Lambda}}$, then it follows by (5.2) that $\mathfrak{b}_N < \infty$ for all $N$. Fix $N \in \mathbb{Z}_+$ and consider the products on the right-hand side of (5.10). Note that

$$\varepsilon \hat{\nu}_a^{(\boldsymbol{\alpha})}(dx_1)(J^{(i_1)}(x_1, dx_2) \cdots J^{(i_n)}(x_{n-1}, dx_n)),$$

consists of a product of blocks which each have the form

$$\varepsilon \hat{\nu}_a^{(\boldsymbol{\alpha})}(dx_j)(\hat{P}_{\boldsymbol{\alpha}}(x_j, dx_{j+1}) \cdots \hat{P}_{\boldsymbol{\alpha}}(x_{j+l-1}, dx_{j+l})) \qquad \text{where } j \ge 1 \text{ and } l \ge 0.$$

(When $l = 0$, it is understood that the product involving the $\hat{P}_{\boldsymbol{\alpha}}$ terms is empty. Also, the leading "$\varepsilon$" term only appears for $j > 1$.) Now the total number of blocks is $k$, where $k$ is the cardinality of $\{j : i_j = 0, 1 \le j < n\} + 1$. Moreover for each block,

$$
\begin{aligned}
(5.13) \qquad &\int_{\mathbb{S}^l} \hat{\nu}_a^{(\boldsymbol{\alpha})}(dx_j)(\hat{P}_{\boldsymbol{\alpha}}(x_j, dx_{j+1}) \cdots \hat{P}_{\boldsymbol{\alpha}}(x_{j+l-1}, dx_{j+l})) \\
&= \mathbb{E}_{\nu_a}^{(1)}[e^{\langle \boldsymbol{\alpha}, \mathbf{S}_{l+1} \rangle}] \le (\sigma_{l+1}(\boldsymbol{\alpha}))^{l+1},
\end{aligned}
$$

where the last step follows from (5.2) and the definition of $\sigma_l$ given just prior to (5.12). Now by definition, $\sigma_l(\boldsymbol{\alpha})$ is decreasing in $l$ for any fixed $\boldsymbol{\alpha}$. Hence if $l + 1 \geq N$ then $\sigma_{l+1}(\boldsymbol{\alpha}) \leq \sigma_N(\boldsymbol{\alpha})$, while if $l + 1 < N$ then $\sigma_{l+1}(\boldsymbol{\alpha}) \leq \sigma_1(\boldsymbol{\alpha})$, and hence

$$(\sigma_{l+1}(\boldsymbol{\alpha}))^{l+1} \leq \mathfrak{b}_N(\sigma_N(\boldsymbol{\alpha}))^{l+1}$$

by (5.12). Substituting these estimates into (5.13) yields

$$(5.14) \quad \int_{\mathbb{S}^{l+1}} \hat{\nu}_a^{(\boldsymbol{\alpha})}(dx_j)(\hat{P}_{\boldsymbol{\alpha}}(x_j, dx_{j+1}) \cdots \hat{P}_{\boldsymbol{\alpha}}(x_{j+l-1}, dx_{j+l})) \leq \mathfrak{b}_N(\sigma_N(\boldsymbol{\alpha}))^{l+1}.$$

If the total number of such blocks attached to a given $(i_1, \ldots, i_{n-1}) \in \mathbb{I}$ is denoted by $k$, then we obtain for each individual integral on the right-hand side of (5.10) that

$$(5.15) \quad \begin{aligned} \int_{\mathbb{S}^n} \hat{\nu}_a^{(\boldsymbol{\alpha})}(dx_1) J^{(i_1)}(x_1, dx_2) \cdots J^{(i_{n-1})}(x_{n-1}, dx_n) \\ \leq \varepsilon^{k-1}(\mathfrak{b}_N\sigma_N(\boldsymbol{\alpha}))^k(\sigma_N(\boldsymbol{\alpha}))^{n-k}. \end{aligned}$$

Summing over all $(i_1, \ldots, i_{n-1}) \in \mathbb{I}$, it follows from (5.10) and (5.15) that

$$(5.16) \quad \mathfrak{L}_n^{(\varepsilon)}(\boldsymbol{\alpha}) \leq \frac{1}{\varepsilon}(\sigma_N(\boldsymbol{\alpha}) + \varepsilon\mathfrak{b}_N\sigma_N(\boldsymbol{\alpha}))^n.$$

Hence

$$(5.17) \quad \limsup_{n \to \infty} \frac{1}{n} \log \mathfrak{L}_n^{(\varepsilon)}(\boldsymbol{\alpha}) \leq \log(\sigma_N(\boldsymbol{\alpha})(1 + \varepsilon\mathfrak{b}_N)).$$

Now let $\varepsilon \to 0$ and then $N \to \infty$. Note that

$$\log \sigma_N(\boldsymbol{\alpha}) = \bar{\boldsymbol{\Lambda}}_N(\boldsymbol{\alpha}) \to \boldsymbol{\Lambda}_{\nu_a}(\boldsymbol{\alpha}) \qquad \text{as } N \to \infty.$$

Consequently (5.11) follows from (5.17).

(iv) This is a variant of Nummelin (1984), Theorem 5.1. First set $\gamma = (\theta^{(\varepsilon)}(\alpha))^{-1}$, and set

$$\mathcal{K} = \gamma(\hat{P}_\alpha^{(\varepsilon)} - \delta_a \mathbf{1}_{\mathcal{L}_a h} \otimes \hat{\nu}_a^{(\alpha)}) \quad \text{and} \quad G = \sum_{n=0}^{\infty} \mathcal{K}^n.$$

Then Proposition 2.1 of Nummelin (1984) gives that $G = I + \mathcal{K}G$, and thus

$$(5.18) \quad \gamma G\delta_a \mathbf{1}_{\mathcal{L}_a h} = \gamma\delta_a \mathbf{1}_{\mathcal{L}_a h} + \gamma^2(\hat{P}_\alpha^{(\varepsilon)} - \delta_a \mathbf{1}_{\mathcal{L}_a h} \otimes \hat{\nu}_a^{(\alpha)})G\delta_a \mathbf{1}_{\mathcal{L}_a h}.$$

The first term on the left-hand side is

$$(5.19) \quad \gamma G\delta_a \mathbf{1}_{\mathcal{L}_a h} := \sum_{n=0}^{\infty} \gamma^{n+1}(\hat{P}_\alpha^{(\varepsilon)} - \delta_a \mathbf{1}_{\mathcal{L}_a h} \otimes \hat{\nu}_a^{(\alpha)})^n \delta_a \mathbf{1}_{\mathcal{L}_a h} := r_\alpha^{(\varepsilon)}.$$

This equivalence also identifies the middle term on the right-hand side of (5.18); namely,

$$(5.20) \qquad \gamma^2 \hat{P}_\alpha^{(\varepsilon)} G \delta_a \mathbf{1}_{\mathcal{L}_a h} = \gamma \hat{P}_\alpha^{(\varepsilon)} r_\alpha^{(\varepsilon)}.$$

Moreover, the last term on the right-hand side of (5.18) is

$$(5.21) \quad \begin{aligned} &\gamma^2 (\delta_a \mathbf{1}_{\mathcal{L}_a h} \otimes \hat{\nu}_a^{(\alpha)}) G \delta_a \mathbf{1}_{\mathcal{L}_a h} \\ &= \gamma \delta_a \mathbf{1}_{\mathcal{L}_a h} \left\{ \sum_{n=1}^\infty \gamma^n \hat{\nu}_a^{(\alpha)} (\hat{P}_\alpha^{(\varepsilon)} - \delta_a \mathbf{1}_{\mathcal{L}_a h} \otimes \hat{\nu}_a^{(\alpha)})^{n-1} \delta_a \mathbf{1}_{\mathcal{L}_a h} \right\} \\ &\begin{cases} = \gamma \delta_a \mathbf{1}_{\mathcal{L}_a h}, & \text{if } \hat{P}_\alpha^{(\varepsilon)} \text{ is } \gamma\text{-recurrent}, \\ < \gamma \delta_a \mathbf{1}_{\mathcal{L}_a h}, & \text{if } \hat{P}_\alpha^{(\varepsilon)} \text{ is } \gamma\text{-transient}, \end{cases} \end{aligned}$$

where the final step was obtained from Nummelin (1984), Proposition 4.7. [In Nummelin's notation, the quantity in brackets is identified as $\hat{b}(\gamma)$, which is equal to one in the $\gamma$-recurrent case and is less than one in the $\gamma$-transient case.] Substituting (5.19), (5.20) and (5.21) into (5.18) yields

$$(5.22) \qquad r_\alpha^{(\varepsilon)} = \gamma \hat{P}_\alpha^{(\varepsilon)} r_\alpha^{(\varepsilon)}, \qquad \text{if } \hat{P}_\alpha^{(\varepsilon)} \text{ is } \gamma\text{-recurrent},$$

and

$$(5.23) \qquad r_\alpha^{(\varepsilon)} > \gamma \hat{P}_\alpha^{(\varepsilon)} r_\alpha^{(\varepsilon)}, \qquad \text{if } \hat{P}_\alpha^{(\varepsilon)} \text{ is } \gamma\text{-transient}.$$

Hence, $r_\alpha^{(\varepsilon)}$ is $\gamma$-invariant when $\hat{P}_\alpha^{(\varepsilon)}$ is $\gamma$-recurrent and $\gamma$-subinvariant when $\hat{P}_\alpha^{(\varepsilon)}$ is $\gamma$-transient. Moreover, after rearranging terms in (5.22) and (5.23), we obtain by definition that

$$Q_\alpha^{(\varepsilon)}(x, \mathbb{S}) := \frac{\gamma (\hat{P}_\alpha^{(\varepsilon)} r_\alpha^{(\varepsilon)})(x)}{r_\alpha^{(\varepsilon)}(x)} \begin{cases} = 1, & \text{if } \hat{P}_\alpha^{(\varepsilon)} \text{ is } \gamma\text{-recurrent}, \\ < 1, & \text{if } \hat{P}_\alpha^{(\varepsilon)} \text{ is } \gamma\text{-transient}. \end{cases}$$

Hence, $Q_\alpha^{(\varepsilon)}$ is a probability kernel in the $\gamma$-recurrent case and a subprobability measure in the $\gamma$-transient case. $\quad\square$

REMARK 5.1. In Proposition 5.1(ii), it has just been observed that

$$\Theta(\alpha) \le \Lambda(\alpha) = 0 \qquad \forall \alpha,$$

but the reverse inequality need not be true, in general. However, under ($\mathfrak{M}$), ($H_2$) and ($H_3$), it can be shown that

$$(5.24) \qquad \Theta(\mathfrak{r}) = 0 = \Lambda(\mathfrak{r}).$$

This last equation can be obtained as a consequence of Theorem 4.2 and the proof of Theorem 2.1.

To establish (5.24) based on the aforementioned results, first note that $\Lambda$ is convex [by Hölder's inequality] and finite in a neighborhood of $\mathfrak{r}$ [by (H$_2$)], since $\boldsymbol{\Lambda}((\alpha, \beta))$ is nondecreasing in $\beta$. Consequently $\Lambda$ is continuous at $\mathfrak{r}$. The definition of $\mathfrak{r}$ then implies that $\Lambda(\mathfrak{r}) = 0$; thus, it is sufficient to show that $\Theta(\mathfrak{r}) = 0$.

To this end, note that in the proof of Theorem 2.1, we have shown that

$$\mathfrak{r} = \eta := \sup\{\alpha : \log \mathbb{E}_{\nu_a}[e^{\alpha \tilde{S}}] \le 0\} \qquad \text{where } \tilde{S} \stackrel{d}{=} S_{T_i - 1} - S_{T_{i-1} - 1} \text{ for } i \ge 1.$$

Let

$$\tilde{\mathfrak{r}} = \sup\{\alpha : \Theta(\alpha) \le 0\}.$$

Now if $\tilde{\mathfrak{r}} > \mathfrak{r}$, then for any $\alpha \in (\mathfrak{r}, \tilde{\mathfrak{r}})$ we would have $\Theta(\alpha) \le 0$ [since $\Theta$ is convex and $\alpha < \tilde{\mathfrak{r}}$], and $\mathbb{E}_{\nu_a}[e^{\alpha \tilde{S}}] > 1$ [since $\alpha > \mathfrak{r}$]. Hence

$$(5.25) \qquad \mathbb{E}_{\nu_a}[e^{\alpha \tilde{S} - \tau \Theta(\alpha)}] \ge \mathbb{E}_{\nu_a}[e^{\alpha \tilde{S}}] > 1$$

and this is a contradiction to (4.24). We conclude that $\tilde{\mathfrak{r}} \le \mathfrak{r}$, which means that $\Theta(\alpha) > 0$ for all $\alpha > \mathfrak{r}$. Now $\Theta \le \Lambda$, and so (H$_2$) $\Longrightarrow \mathfrak{r} \in \mathrm{int}\, \mathfrak{D}_\Theta$. Since a convex function is continuous on the interior of its domain, it follows that $\Theta$ is continuous at $\mathfrak{r}$. Moreover, as we have observed, $\Theta(\mathfrak{r}) \le 0$, and $\Theta(\alpha) > 0$ for all $\alpha > \mathfrak{r}$. By continuity, we conclude that $\Theta(\mathfrak{r}) = 0$, thus establishing (5.24). $\qquad \square$

PROPOSITION 5.2. *Suppose that* ($\mathcal{M}$) *and* (H$_3$) *are satisfied, and suppose that there exist $\beta > 0$ and $\varepsilon > 0$ such that $(\mathfrak{r}, \beta) \in \mathrm{int}\, \mathfrak{D}_{\boldsymbol{\Lambda}^{(\varepsilon)}}$. Then:*

(i) *The functions $r_\alpha^{(\varepsilon)}$ and $\mathbf{r}_{\boldsymbol{\alpha}}^{(\varepsilon)}$ are uniformly bounded from above on $\mathcal{L}_a h$. Moreover, there exists a sequence $\{(\alpha_j, \varepsilon_j)\}_{j \in \mathbb{Z}_+}$, where $(\alpha_j, \varepsilon_j) \searrow (\mathfrak{r}, 0)$ as $j \to \infty$, and the sequence $\{\boldsymbol{\alpha}_j\} := \{(\alpha_j, \beta)\}$, such that if*

$$\mathcal{C}_j := \sup_{x \in \mathcal{L}_a h} r_{\alpha_j}^{(\varepsilon_j)}(x) \quad \text{and} \quad \mathcal{D}_j := \sup_{x \in \mathcal{L}_a h} \mathbf{r}_{\boldsymbol{\alpha}_j}^{(\varepsilon_j)}(x) \qquad \forall j,$$

*then $\{\mathcal{C}_j\}_{j \in \mathbb{Z}_+}$ and $\{\mathcal{D}_j\}_{j \in \mathbb{Z}_+}$ are bounded from above by finite constants.*

(ii) *The functions $r_\alpha^{(\varepsilon)}$ and $\mathbf{r}_{\boldsymbol{\alpha}}^{(\varepsilon)}$ are uniformly positive on $\mathbb{S}$. Furthermore, if $\{(\alpha_j, \varepsilon_j)\}$ and $\{\boldsymbol{\alpha}_j\}$ are given as in* (i) *and*

$$\tilde{\mathcal{C}}_j := \inf_{x \in \mathcal{L}_a h} r_{\alpha_j}^{(\varepsilon_j)}(x) \quad \text{and} \quad \tilde{\mathcal{D}}_j := \inf_{x \in \mathcal{L}_a h} \mathbf{r}_{\boldsymbol{\alpha}_j}^{(\varepsilon_j)}(x) \qquad \forall j,$$

*then, possibly after passing to a subsequence, we have that $\{\tilde{\mathcal{C}}_j\}_{j \in \mathbb{Z}_+}$ and $\{\tilde{\mathcal{D}}_j\}_{j \in \mathbb{Z}_+}$ are bounded from below by positive constants.*

PROOF. (i) Without loss of generality, we may assume that the constant $D_a$ in (H$_3$) is greater than or equal to one.

For any $x \in \mathcal{L}_a h$ and $E \in \mathcal{S}$, it follows from (H$_3$) that

$$(5.26) \quad \hat{P}_\alpha^{(\varepsilon)}(x, E) = \int_E e^{\alpha f(y)} P(x, dy) + \varepsilon \hat{\nu}_a^{(\alpha)}(E) \leq D_a \sum_{i=1}^l \hat{P}_\alpha^{(\varepsilon)}(x_i, E),$$

where $x_i \in E_i$ for all $i$. Consequently, for any $x \in \mathcal{L}_a h$,

$$(5.27) \quad \begin{aligned} \mathcal{K}(x, E) &:= (\hat{P}_\alpha^{(\varepsilon)} - \delta_a \mathbf{1}_{\mathcal{L}_a h} \otimes \hat{\nu}_a^{(\alpha)})(x, E) \\ &\leq D_a \sum_{i=1}^l \mathcal{K}(x_i, E) + (D_a l - 1)\delta_a \hat{\nu}_a^{(\alpha)}(E). \end{aligned}$$

Hence for all $x \in \mathcal{L}_a h$ and $n \geq 1$,

$$(5.28) \quad \begin{aligned} \mathcal{K}^n(x, E) &\leq \int_y \left( D_a \sum_{i=1}^l \mathcal{K}(x_i, dy) + (D_a l - 1)\delta_a \hat{\nu}_a^{(\alpha)}(dy) \right) \mathcal{K}^{n-1}(y, E) \\ &= D_a \sum_{i=1}^l \mathcal{K}^n(x_i, E) + (D_a l - 1)\delta_a(\hat{\nu}_a^{(\alpha)} \mathcal{K}^{n-1})(E). \end{aligned}$$

It follows from (5.28) and the definition of $r_\alpha^{(\varepsilon)}$ that, for any $x \in \mathcal{L}_a h$,

$$(5.29) \quad \begin{aligned} r_\alpha^{(\varepsilon)}(x) &:= \sum_{n=0}^\infty (\theta^{(\varepsilon)}(\alpha))^{-n-1} (\mathcal{K}^n \delta_a \mathbf{1}_{\mathcal{L}_a h})(x) \\ &\leq D_a \sum_{i=1}^l r_\alpha^{(\varepsilon)}(x_i) + \frac{D_a l - 1}{\theta^{(\varepsilon)}(\alpha)}, \end{aligned}$$

where, in the last step, the second quantity on the right-hand side was obtained from Nummelin (1984), Proposition 4.7(ii). Thus $r_\alpha^{(\varepsilon)}$ is bounded from above on $\mathcal{L}_a h$, provided that the quantity on the right-hand side is finite. But $r_\alpha^{(\varepsilon)}$ is necessarily finite $\varphi$-a.s. [Nummelin (1984), Proposition 5.1]. Since the $E_i$-sets in (H$_3$) are $\varphi$-positive, we may then choose $x_1, \ldots, x_l$ such that $r_\alpha^{(\varepsilon)}(x_i)$ is finite for each $i$.

To establish uniform boundedness along a sequence $\{(\alpha_j, \varepsilon_j)\}$, where $(\alpha_j, \varepsilon_j) \searrow (\mathfrak{r}, 0)$ as $j \to \infty$, set

$$(5.30) \quad q_1(x) = \liminf_{(\alpha, \varepsilon) \searrow (\mathfrak{r}, 0)} r_\alpha^{(\varepsilon)}(x) \qquad \forall x \in \mathbb{S},$$

where the limit is taken along the straight-line path joining $(\mathfrak{r}, 0)$ to some point $(\hat{\alpha}, \hat{\varepsilon})$ with $\hat{\alpha} > \mathfrak{r}$ and $\hat{\varepsilon} > 0$.

We begin by identifying the limit of $\Theta^{(\varepsilon)}(\alpha)$ as $(\alpha, \varepsilon) \searrow (\mathfrak{r}, 0)$. It follows immediately from the definition of $\Theta^{(\varepsilon)}$ that $\Theta^{(\varepsilon)}(\alpha) \geq \Theta(\alpha)$ for all $(\alpha, \varepsilon)$; and hence

$$(5.31) \qquad \liminf_{(\alpha, \varepsilon) \searrow (\mathfrak{r}, 0)} \Theta^{(\varepsilon)}(\alpha) \geq \lim_{\alpha \searrow \mathfrak{r}} \Theta(\alpha) = \Theta(\mathfrak{r}) > -\infty \qquad \forall \alpha, \; \forall \varepsilon > 0,$$

where the last two inequalities follow from Theorem 3.2 of Nummelin (1984) and the convexity of $\Theta$. Moreover, it follows by Proposition 5.1 that

$$(5.32) \qquad \Theta^{(\varepsilon)}(\alpha) \leq \Lambda^{(\varepsilon)}(\alpha) \searrow \Lambda(\mathfrak{r}) \qquad \text{as } (\alpha, \varepsilon) \searrow (\mathfrak{r}, 0);$$

and from the definition of $\mathfrak{r}$, we also have that $\Lambda(\mathfrak{r}) \leq 0$. [In the last step of (5.32), we have used the fact that $\Lambda$ is convex, which follows from an application of Hölder's inequality, and we have used our assumption that $(\mathfrak{r}, 0) \in \operatorname{int} \mathfrak{D}_{\Lambda}$, which implies that $\mathfrak{r} \in \operatorname{int} \mathfrak{D}_{\Lambda}$.] From (5.31) and (5.32) we conclude that

$$(5.33) \qquad \liminf_{(\alpha, \varepsilon) \searrow (\mathfrak{r}, 0)} \theta^{(\varepsilon)}(\alpha) = \mathfrak{a} \in (\theta(\mathfrak{r}), 1].$$

Recalling the definition of $q_1$ in (5.30), now apply Fatou's lemma to the equation

$$\int_{\mathbb{S}} e^{\alpha f(y)} r_\alpha^{(\varepsilon)}(y) P^{(\varepsilon)}(x, dy) \leq \theta^{(\varepsilon)}(\alpha) r_\alpha^{(\varepsilon)}(x)$$

to obtain

$$(5.34) \qquad \int_{\mathbb{S}} e^{\mathfrak{r} f(y)} q_1(y) P(x, dy) \leq q_1(x).$$

Thus $q_1$ is a subinvariant function for the kernel $\hat{P}_{\mathfrak{r}}$. Also, by Nummelin (1984), Proposition 5.1, the set $\{x : q_1(x) < \infty\}$ is full; hence, either $q_1 < \infty$ $\varphi$-a.s. or $q_1 \equiv \infty$.

To see that the latter is not the case, observe by the definition of $r_\alpha^{(\varepsilon)}$ and Proposition 4.7(ii) of Nummelin (1984) that

$$\int_{\mathbb{S}} r_\alpha^{(\varepsilon)}(x) \hat{\nu}_a^{(\alpha)}(dx) \leq 1.$$

Applying Fatou's lemma to this equation yields

$$(5.35) \qquad \int_{\mathbb{S}} q_1(y) \hat{\nu}_a^{(\mathfrak{r})}(dx) \leq 1.$$

Consequently $q_1 \not\equiv \infty$.

Hence, if $E_1$ is given as in (H3), then

$$(5.36) \qquad \liminf_{(\alpha, \varepsilon) \searrow (\mathfrak{r}, 0)} r_\alpha^{(\varepsilon)}(x_1) := q_1(x_1) < \infty, \qquad \varphi\text{-a.a. } x_1 \in E_1.$$

Therefore, for any given $x_1 \in E_1$, there exists a subsequence $\{(\alpha'_j, \varepsilon'_j)\}_{j \in \mathbb{Z}_+}$ such that

$$(5.37) \qquad \sup_j r^{(\varepsilon'_j)}_{\alpha'_j}(x_1) < \infty.$$

Next, repeat the same argument but with

$$q_2(x) := \liminf_{j \to \infty} r^{(\varepsilon'_j)}_{\alpha'_j}(x)$$

in place of $q_1$, to establish the existence of an element $x_2 \in E_2$ and a subsequence $\{(\alpha''_j, \varepsilon''_j)\} \subseteq \{(\alpha'_j, \varepsilon'_j)\}$ such that

$$(5.38) \qquad \sup_j r^{(\varepsilon''_j)}_{\alpha''_j}(x_2) < \infty.$$

Continuing in this manner, we obtain elements $x_i \in E_i$, $i = 1, \ldots, l$, and a sequence $\{(\alpha_j, \varepsilon_j)\}$ where $(\alpha_j, \varepsilon_j) \searrow (\mathfrak{r}, 0)$ as $j \to \infty$, such that

$$(5.39) \qquad \sup_j r^{(\varepsilon_j)}_{\alpha_j}(x_i) < \infty, \qquad i = 1, \ldots, l.$$

Substituting this expression and (5.33) into (5.29) yields that for some positive integer $J_0$,

$$\sup_{j \geq J_0} r^{(\varepsilon_j)}_{\alpha_j}(x) \leq D_a l \left( \sup_{i,j} r^{(\varepsilon_j)}_{\alpha_j}(x_i) + \frac{1}{\mathfrak{a}} \right) \qquad \forall x \in \mathcal{L}_a h.$$

Since the quantity on the right-hand side is independent of $x \in \mathcal{L}_a h$ and $j \geq J_0$, we conclude that $\{\mathcal{C}_j\}_{j \geq J_0}$ is bounded.

The proof of these properties for the function $\mathbf{r}^{(\varepsilon)}_{\boldsymbol{\alpha}}$ is essentially the same. The main modification arises in the definition of $q_1$, which is now replaced with

$$\mathbf{q}_1(x) := \liminf_{j \to \infty} \mathbf{r}^{(\varepsilon_j)}_{\boldsymbol{\alpha}_j}(x),$$

where $\{(\alpha_j, \varepsilon_j)\}$ is given as above. The previous argument may then be repeated to obtain successive subsequences corresponding to functions $\mathbf{q}_1$, $\mathbf{q}_2, \ldots$ . (The choice of the elements $x_i \in E_i$ may, of course, be different in this case.) As before, we obtain boundedness for the sequence $\{\mathcal{D}_j\}$, but now along a subsequence $\{(\alpha_{j_k}, \varepsilon_{j_k})\}$ of $\{(\alpha_j, \varepsilon_j)\}$. Thus (i) holds.

(ii) To show that $r^{(\varepsilon)}_\alpha$ is uniformly positive, note by definition that

$$\hat{P}^{(\varepsilon)}(x, E) \geq \varepsilon \int_E e^{\alpha f(y)} \nu_a(dy) := \varepsilon \hat{\nu}^{(\alpha)}_a(E).$$

Since $r_\alpha^{(\varepsilon)}$ is $(\theta^{(\varepsilon)}(\alpha))^{-1}$-subinvariant with respect to the kernel $\hat{P}_\alpha^{(\varepsilon)}$, it follows that

$$
\begin{aligned}
(5.40) \qquad r_\alpha^{(\varepsilon)}(x) &\geq (\theta^{(\varepsilon)}(\alpha))^{-1} \int_{\mathbb{S}} r_\alpha^{(\varepsilon)}(y) \hat{P}_\alpha^{(\varepsilon)}(x, dy) \\
&\geq (\theta^{(\varepsilon)}(\alpha))^{-1} \varepsilon \int_{\mathbb{S}} r_\alpha^{(\varepsilon)}(y) \hat{\nu}_a^{(\alpha)}(dy) = \text{const.} > 0,
\end{aligned}
$$

where the last inequality was obtained from the positivity of $r_\alpha^{(\varepsilon)}$. Hence $r_\alpha^{(\varepsilon)}$ is uniformly positive.

To establish the uniform positivity of $\{\tilde{\mathcal{C}}_j\}$, note by the definition of $r_\alpha^{(\varepsilon)}$ and the definition of $\mathcal{K}$ [in (5.27)] that

$$
(5.41) \qquad r_\alpha^{(\varepsilon)}(x) := \sum_{n=0}^{\infty} (\theta^{(\varepsilon)}(\alpha))^{-n-1} (\mathcal{K}^n \delta_a \mathbf{1}_{\mathcal{L}_a h})(x) \geq \frac{\delta_a \mathbf{1}_{\mathcal{L}_a h}(x)}{\theta^{(\varepsilon)}(\alpha)}.
$$

Since $\liminf_{j \to \infty} \theta^{(\varepsilon_j)}(\alpha_j) = \mathfrak{a} \in (\theta(\mathfrak{r}), 1]$ [by (5.33)], it follows that for some positive integer $J_1 \geq J_0$,

$$
(5.42) \qquad \inf_{j \geq J_1} r_{\alpha_j}^{(\varepsilon_j)}(x) \geq \frac{\delta_a}{2\mathfrak{a}} > 0 \qquad \forall x \in \mathcal{L}_a h.
$$

Hence—possibly after redefining the sequence $\{(\alpha_j, \varepsilon_j)\}$ so that $(\alpha_{J_1}, \varepsilon_{J_1})$ is actually the initial term of this sequence—we obtain that $\{\tilde{\mathcal{C}}_j\}$ bounded from below by a positive constant.

The same reasoning shows that the sequence $\{\tilde{\mathcal{D}}_j\}$ is likewise bounded from below by a positive constant.  $\square$

5.3. *Proof of Theorem 4.2.*  The proof will be based on two lemmas. To motivate the first of these lemmas, recall that a Markov chain is *geometric recurrent* if

$$
(\mathcal{G}) \qquad\qquad \sup_{x \in \mathcal{L}_a h} \mathbb{E}_x[\rho^\tau] < \infty \qquad \text{for some } \rho > 1,
$$

where $\tau$ is a typical regeneration time. To establish $(\mathcal{G})$, one usually begins by showing that $(\mathcal{G})$ holds under the additional assumption that there exists an atom, namely,

$$
(\mathcal{G}') \qquad\qquad \sup_{x \in \mathcal{L}_a h} \mathbb{E}_x[\rho^\mathfrak{T}] < \infty \qquad \text{for some } \rho > 1,
$$

where $\mathfrak{T}$ denotes the first return time of the chain to its regeneration set, that is, $\mathfrak{T} := \inf\{n : X_n \in \mathcal{L}_a h\}$.

Let $Q_\alpha, r_\alpha$ be defined as $Q_\alpha^{(\varepsilon)}, r_\alpha^{(\varepsilon)}$ but with $\varepsilon = 0$; cf. Section 5.1. Then the goal of our first lemma will be to establish $(\mathcal{G}')$ for the *shifted* kernel $Q_\alpha$; that is,

$$
(5.43) \qquad \sup_{x \in \mathcal{L}_a h} \mathbb{E}_x^Q[\rho^\mathfrak{T}] := \sup_{x \in \mathcal{L}_a h} \mathbb{E}_x\left[ \frac{\rho^\mathfrak{T} e^{\alpha S_\mathfrak{T} - \mathfrak{T}\Theta(\alpha)} r_\alpha(X_\mathfrak{T})}{r_\alpha(X_0)} \right] < \infty,
$$

where $\mathbb{E}^Q[\cdot]$ denotes the expectation under the measure $Q_\alpha$ (and for this heuristical discussion, we assume that $Q_\alpha$ is indeed a probability measure). Since $\Theta(\alpha) \searrow \Theta(\mathfrak{r}) \approx 0$ as $\alpha \to \mathfrak{r}$, it is a natural strengthening of $(\mathcal{G}')$ to require that

$$(5.44) \qquad \sup_{x \in \mathcal{L}_a h} \mathbb{E}_x \left[ \frac{\rho^{\mathfrak{T}} e^{\alpha S_{\mathfrak{T}}} r_\alpha(X_{\mathfrak{T}})}{r_\alpha(X_0)} \right] < \infty \qquad \text{for some } \alpha > \mathfrak{r}.$$

While (5.44) has a straightforward probabilistic interpretation, it is not sufficiently strong for our purposes because we will ultimately need the eigenfunctions $r_\alpha$ to be uniformly positive, which may not be true in general. Hence we will need to work with the perturbed kernel $P^{(\varepsilon)}$ rather than $P$. But $P^{(\varepsilon)}$ is not a probability kernel, and consequently, we must first replace (5.44) with a more general series representation. To this end let

$$(5.45) \qquad \mathfrak{P}^{(\varepsilon)}(x, E) = \int_{E \cap (\mathcal{L}_a h)^c} P^{(\varepsilon)}(x, dy)$$

and

$$(5.46) \qquad \underline{P}^{(\varepsilon)}(x, E) = \int_{E \cap \mathcal{L}_a h} P^{(\varepsilon)}(x, dy) \qquad \forall \varepsilon \geq 0,$$

and set $\mathfrak{P} = \mathfrak{P}^{(0)}$ and $\underline{P} = \underline{P}^{(0)}$. Qualitatively, $\mathfrak{P}$ is a kernel which is nonzero and equal to $P$ when $\{X_n\}$ *avoids* the regeneration set $\mathcal{L}_a h$, and $\underline{P}$ is a kernel which is nonzero when $\{X_n\}$ *enters* $\mathcal{L}_a h$. Thus, in particular, $\mathfrak{P}^{n-1} \underline{P}(x, \mathbb{S})$ describes the event that $\{X_n\}$ begins at state $x$, avoids the set $\mathcal{L}_a h$ during its first $n-1$ time steps, and enters the set $\mathcal{L}_a h$ during the $n$th time step. Thus it follows from the definitions that

$$(5.47) \qquad \mathbb{E}_x \left[ \frac{\rho^{\mathfrak{T}} e^{\alpha S_{\mathfrak{T}}} r_\alpha(X_{\mathfrak{T}})}{r_\alpha(X_0)} \right] = \frac{1}{r_\alpha(x)} \sum_{n=1}^{\infty} \rho^n (\hat{\mathfrak{P}}_\alpha^{n-1} \hat{\underline{P}}_\alpha r_\alpha)(x).$$

Our first objective is to establish the following variant of (5.44), (5.47).

LEMMA 5.1. *Assume that* $(\mathfrak{M})$, $(\mathrm{H}_2)$ *and* $(\mathrm{H}_3)$ *are satisfied. Then for sufficiently large $a > 0$ and any constant $\rho > 0$, there exist a constant $D < \infty$ such that*

$$(5.48) \qquad \sup_{x \in \mathcal{L}_a h} \left\{ \frac{1}{r_\alpha^{(\varepsilon)}(x)} \sum_{n=1}^{\infty} \rho^n ((\hat{\mathfrak{P}}_\alpha^{(\varepsilon)})^{n-1} \hat{\underline{P}}_\alpha^{(\varepsilon)} r_\alpha^{(\varepsilon)})(x) \right\} \leq D,$$

*uniformly for* $(\alpha, \varepsilon) \in \{(\alpha_j, \varepsilon_j)\}_{j \geq J}$, *where* $(\alpha_j, \varepsilon_j) \searrow (\mathfrak{r}, 0)$ *is given as in Proposition 5.2 and $J$ is a positive integer.*

PROOF. By $(\mathrm{H}_2)$, there exist points $\hat{\alpha} > \mathfrak{r}$ and $\hat{\beta} > 0$ such that $\mathbf{\Lambda}(\hat{\alpha}, \hat{\beta}) < \infty$. Set $\hat{\boldsymbol{\alpha}} = (\hat{\alpha}, \hat{\beta})$, and observe by Proposition 5.1 that

$$\mathbf{\Theta}^{(\varepsilon)}(\hat{\boldsymbol{\alpha}}) \leq \mathbf{\Lambda}^{(\varepsilon)}(\hat{\boldsymbol{\alpha}}) \searrow \mathbf{\Lambda}(\hat{\boldsymbol{\alpha}}) \qquad \text{as } \varepsilon \to 0.$$

Hence

$$\boldsymbol{\Theta}^{(\varepsilon)}(\hat{\boldsymbol{\alpha}}) < \infty \qquad \text{for } \varepsilon \geq 0 \text{ sufficiently small.}$$

By a similar argument,

$$\boldsymbol{\Theta}^{(\varepsilon)}(\mathbf{0}) \leq 1 \qquad \text{for } \varepsilon \geq 0 \text{ sufficiently small.}$$

Then it follows from the convexity of $\boldsymbol{\Theta}^{(\varepsilon)}$ that on the line segment joining the origin to $(\hat{\alpha}, \hat{\beta})$, the function $\boldsymbol{\Theta}^{(\varepsilon)}$ always lies below $\max\{\boldsymbol{\Theta}^{(\varepsilon)}(\mathbf{0}), \boldsymbol{\Theta}^{(\varepsilon)}(\hat{\alpha}, \hat{\beta})\}$. In other words, for some $\hat{\varepsilon} > 0$ we have

$$(5.49) \qquad \boldsymbol{\Theta}^{(\hat{\varepsilon})}\left(\alpha, \frac{\hat{\beta}}{\hat{\alpha}}\alpha\right) \leq \max\{1, \boldsymbol{\Theta}^{(\hat{\varepsilon})}(\hat{\alpha}, \hat{\beta})\} < \infty \qquad \forall \alpha \in [0, \hat{\alpha}].$$

Now by definition, $\boldsymbol{\Theta}^{(\varepsilon)}(\alpha, \beta)$ is nondecreasing in $\beta$ [for fixed $(\alpha, \varepsilon)$ and increasing $\beta$], and $\boldsymbol{\Theta}^{(\varepsilon)}(\alpha, \beta)$ is nondecreasing in $\varepsilon$ [for fixed $(\alpha, \beta)$ and increasing $\varepsilon$]. Consequently it follows from (5.49) that

$$(5.50) \qquad \boldsymbol{\Theta}^{(\varepsilon)}(\alpha, \beta) \leq \max\{1, \boldsymbol{\Theta}^{(\hat{\varepsilon})}(\hat{\alpha}, \hat{\beta})\} < \infty$$

$$\forall \alpha \in [\mathfrak{r}, \hat{\alpha}], \ \varepsilon \in [0, \hat{\varepsilon}], \text{ and for } \beta = \frac{\hat{\beta}}{\hat{\alpha}}\mathfrak{r}.$$

In other words, $\boldsymbol{\Theta}^{(\varepsilon)}(\alpha, \beta)$ is bounded *uniformly* when $(\alpha, \varepsilon) \in [\mathfrak{r}, \hat{\alpha}] \times [0, \hat{\varepsilon}]$. Therefore, for any $\rho > 0$, and for $\boldsymbol{\alpha} = (\alpha, \beta)$ and $\beta$ is given as in (5.50), there exists a finite constant $a$ such that

$$(5.51) \qquad \boldsymbol{\theta}^{(\varepsilon)}(\boldsymbol{\alpha}) := \exp \boldsymbol{\Theta}^{(\varepsilon)}(\boldsymbol{\alpha}) \leq \frac{e^{\beta a}}{2\rho} \qquad \forall (\alpha, \varepsilon) \in [\mathfrak{r}, \hat{\alpha}] \times [0, \hat{\varepsilon}].$$

For the remainder of the proof, the constants $\beta$ and $a$ will now be *fixed*, while the parameters $\alpha$ and $\varepsilon$ will be allowed to vary within the range of values specified in (5.51). [Later, we will restrict $(\alpha, \varepsilon)$ to those elements belonging to the sequence $\{(\alpha_j, \varepsilon_j)\}$ which appeared in the statement of Proposition 5.2.]

By the subinvariance of $\mathbf{r}_{\boldsymbol{\alpha}}^{(\varepsilon)}$,

$$(5.52) \qquad \int_{\mathbb{S}} e^{\alpha f(y) + \beta h(y)} P^{(\varepsilon)}(x, dy) \mathbf{r}_{\boldsymbol{\alpha}}^{(\varepsilon)}(y) \leq \boldsymbol{\theta}^{(\varepsilon)}(\boldsymbol{\alpha}) \mathbf{r}_{\boldsymbol{\alpha}}^{(\varepsilon)}(x).$$

Hence

$$(5.53) \qquad e^{\beta a} \int_{(\mathcal{L}_a h)^c} e^{\alpha f(y)} P^{(\varepsilon)}(x, dy) \mathbf{r}_{\boldsymbol{\alpha}}^{(\varepsilon)}(y) \leq \boldsymbol{\theta}^{(\varepsilon)}(\boldsymbol{\alpha}) \mathbf{r}_{\boldsymbol{\alpha}}^{(\varepsilon)}(x).$$

By setting

$$\mathcal{K}(x, E) = \int_{E \cap (\mathcal{L}_a h)^c} \rho \tilde{P}_{\alpha}^{(\varepsilon)}(x, dy) \quad \text{and} \quad V = \mathbf{r}_{\boldsymbol{\alpha}}^{(\varepsilon)}$$

and substituting (5.51) into (5.53), we then obtain

$$\mathcal{K}V \leq \frac{1}{2}V. \tag{5.54}$$

Now if $(\alpha, \varepsilon) \in \{(\alpha_j, \varepsilon_j)\}$, then it follows by Proposition 5.2 that $r_\alpha^{(\varepsilon)}$ is uniformly bounded from below on $\mathcal{L}_a h$. Hence for some positive constant $c$,

$$
\begin{aligned}
\boldsymbol{\theta}^{(\varepsilon)}(\boldsymbol{\alpha}) \mathbf{r}_{\boldsymbol{\alpha}}^{(\varepsilon)}(x) &\geq \int_{\mathcal{L}_a h} \hat{P}_{\boldsymbol{\alpha}}^{(\varepsilon)}(x, dy) \mathbf{r}_{\boldsymbol{\alpha}}^{(\varepsilon)}(y) \geq c \int_{\mathcal{L}_a h} \hat{P}_{\boldsymbol{\alpha}}^{(\varepsilon)}(x, dy) \\
&\geq c \int_{\mathbb{S}} \underline{\hat{P}}_{\alpha}^{(\varepsilon)}(x, dy),
\end{aligned} \tag{5.55}
$$

where, in the last inequality, we have used the fact that $h \geq 0 \Longrightarrow \hat{P}_{\boldsymbol{\beta}} \geq \hat{P}_\alpha$ for all $\beta \geq 0$ and $\boldsymbol{\alpha} = (\alpha, \beta)$, and we have used the fact that $\underline{\hat{P}}_{\alpha}^{(\varepsilon)}(x, \cdot) = 0$ on $(\mathcal{L}_a h)^c$ [cf. (5.46)]. If we now define the constant $\Delta > 0$ and function $\boldsymbol{\Delta} : \mathbb{S} \to \mathbb{R}$ by

$$\Delta = \frac{c}{2\boldsymbol{\theta}^{(\varepsilon)}(\boldsymbol{\alpha})} \quad \text{and} \quad \boldsymbol{\Delta}(x) = \Delta, \ \forall x \in \mathbb{S}, \tag{5.56}$$

then from (5.55) we obtain

$$\underline{\hat{P}}_{\alpha}^{(\varepsilon)} \boldsymbol{\Delta} \leq \frac{1}{2} \mathbf{r}_{\boldsymbol{\alpha}}^{(\varepsilon)} := \frac{1}{2}V. \tag{5.57}$$

From (5.54) and (5.57), we conclude

$$\mathcal{K}V \leq V - \underline{\hat{P}}_{\alpha}^{(\varepsilon)} \boldsymbol{\Delta}, \tag{5.58}$$

which may roughly be viewed as a drift condition satisfied by the kernel $\mathcal{K}$.

Next, we claim that

$$V - \mathcal{K}^m V \geq \sum_{n=1}^{m} \mathcal{K}^{n-1} \underline{\hat{P}}_{\alpha}^{(\varepsilon)} \boldsymbol{\Delta} \qquad \forall m \in \mathbb{Z}_+. \tag{5.59}$$

To establish this claim, proceed by induction, observing that (5.59) holds when $m = 1$, by (5.58). Now assume that (5.59) holds for arbitrary $m$. Multiply (5.59) on the left by $\mathcal{K}$ and apply (5.58) to the first term on the left-hand side to obtain

$$(V - \underline{\hat{P}}_{\alpha}^{(\varepsilon)} \boldsymbol{\Delta}) - \mathcal{K}^{m+1} V \geq \sum_{n=2}^{m+1} \mathcal{K}^{n-1} \underline{\hat{P}}_{\alpha}^{(\varepsilon)} \boldsymbol{\Delta}, \tag{5.60}$$

which is (5.59) with $m + 1$ in place of $m$.

Since $\mathcal{K}^m V \geq 0$, it follows upon letting $m \to \infty$ in (5.59) that

$$
\begin{aligned}
V(x) &\geq \sum_{n=1}^{\infty} (\mathcal{K}^{n-1} \underline{\hat{P}}_\alpha^{(\varepsilon)} \boldsymbol{\Delta})(x) \\
&:= \sum_{n=1}^{\infty} \rho^{n-1} ((\hat{\mathfrak{P}}_\alpha^{(\varepsilon)})^{n-1} \underline{\hat{P}}_\alpha^{(\varepsilon)} \boldsymbol{\Delta})(x).
\end{aligned}
$$

(5.61)

Now by the definition of $\underline{P}^{(\varepsilon)}$, we have that $\underline{\hat{P}}_\alpha^{(\varepsilon)}(x, \cdot) = 0$ on $(\mathcal{L}_a h)^c$, and by Proposition 5.2, $r_\alpha^{(\varepsilon)}$ is bounded from above on $\mathcal{L}_a h$ for all $(\alpha, \varepsilon) \in \{(\alpha_j, \varepsilon_j)\}$. Hence for any $(\alpha, \varepsilon) \in \{(\alpha_j, \varepsilon_j)\}$, there exists a finite constant $d$ such that

$$
\int_E \underline{\hat{P}}_\alpha^{(\varepsilon)}(x, dy) r_\alpha^{(\varepsilon)}(y) \leq d \int_E \underline{\hat{P}}_\alpha^{(\varepsilon)}(x, dy) = \frac{d}{\Delta} \underline{\hat{P}}_\alpha^{(\varepsilon)} \boldsymbol{\Delta}(x) \qquad \forall E \in \mathcal{S};
$$

that is, $\underline{\hat{P}}_\alpha^{(\varepsilon)} \boldsymbol{\Delta} \geq (\Delta/d) \underline{\hat{P}}_\alpha^{(\varepsilon)} r_\alpha^{(\varepsilon)}$. Substituting this inequality into (5.61) yields

(5.62) $$ \frac{1}{r_\alpha^{(\varepsilon)}(x)} \sum_{n=1}^{\infty} \rho^n ((\hat{\mathfrak{P}}_\alpha^{(\varepsilon)})^{n-1} \underline{\hat{P}}_\alpha^{(\varepsilon)} r_\alpha^{(\varepsilon)})(x) \leq \frac{\rho d}{\Delta r_\alpha^{(\varepsilon)}(x)} V(x). $$

Thus we have established (5.48), except that the constant $D$ on the right-hand side of that equation must be replaced with the quantity on the right-hand side of (5.62), which [from the definition of $\Delta$ given in (5.56)] may be identified as

$$
D_{\boldsymbol{\alpha}}^{(\varepsilon)} = \frac{2\rho d \boldsymbol{\theta}^{(\varepsilon)}(\boldsymbol{\alpha})}{c r_\alpha^{(\varepsilon)}(x)} V(x).
$$

It remains to show that $D_{\boldsymbol{\alpha}}^{(\varepsilon)}$ is uniformly bounded from above when $x \in \mathcal{L}_a h$ and $(\alpha, \varepsilon) \in \{(\alpha_j, \varepsilon_j)\}$. But by (5.51), $\boldsymbol{\theta}^{(\varepsilon)}(\boldsymbol{\alpha})$ is uniformly bounded from above for $\alpha \in [\mathfrak{r}, \hat{\alpha}]$ and $\varepsilon \in [0, \hat{\varepsilon}]$. Moreover, by Proposition 5.2, $r_\alpha^{(\varepsilon)}$ and $V := \mathbf{r}_{\boldsymbol{\alpha}}^{(\varepsilon)}$ are uniformly bounded from above and below on the set $\mathcal{L}_a h$, provided that $(\alpha, \varepsilon) \in \{(\alpha_j, \varepsilon_j)\}$. Thus we conclude that $D_{\boldsymbol{\alpha}}^{(\varepsilon)}$ is uniformly bounded from above for $x \in \mathcal{L}_a h$ and $\alpha \in \{(\alpha_j, \varepsilon_j)\} \cap ([\mathfrak{r}, \hat{\alpha}] \times [0, \hat{\varepsilon}]) = \{(\alpha_j, \varepsilon_j)\}_{j \geq J}$, for some positive integer $J$. $\quad \square$

LEMMA 5.2. *Assume that* (𝔐), (H₂) *and* (H₃) *are satisfied, and let* $S_n' = S_n + f(X_0)$. *Then for* $a > 0$ *sufficiently large and* $\rho > 1$ *sufficiently small,*

(5.63) $$ \mathbb{E}_{\nu_a} \left[ \sum_{n=0}^{\tau-1} \rho^n e^{\alpha S_n'} \right] < \infty \qquad \text{for some } \alpha > \mathfrak{r}. $$

PROOF. For any $x \in \mathbb{S}$ and $E \in \mathcal{S}$, set

$$\mathfrak{Q}_\alpha^{(\varepsilon)}(x, E) = \int_{E \cap (\mathcal{L}_a h)^c} Q^{(\varepsilon)}(x, dy).$$

Then it follows directly from the definitions in Section 5.1 that

$$Q_\alpha^{(\varepsilon)}(x, E) = \int_E \frac{r_\alpha^{(\varepsilon)}(y)}{\theta^{(\varepsilon)}(\alpha) r_\alpha^{(\varepsilon)}(x)} \hat{P}_\alpha^{(\varepsilon)}(x, dy)$$

and

$$\mathfrak{Q}_\alpha^{(\varepsilon)}(x, E) = \int_E \frac{r_\alpha^{(\varepsilon)}(y)}{\theta^{(\varepsilon)}(\alpha) r_\alpha^{(\varepsilon)}(x)} \hat{\mathfrak{P}}_\alpha^{(\varepsilon)}(x, dy).$$

Now by Proposition 5.1, $Q_\alpha^{(\varepsilon)}$ is a probability measure when $\hat{P}_\alpha^{(\varepsilon)}$ is $(\theta^{(\varepsilon)}(\alpha))^{-1}$-recurrent, and $Q_\alpha^{(\varepsilon)}$ is a subprobability measure when $\hat{P}_\alpha^{(\varepsilon)}$ is $(\theta^{(\varepsilon)}(\alpha))^{-1}$-transient. First assume that $\hat{P}_\alpha^{(\varepsilon)}$ is $(\theta^{(\varepsilon)}(\alpha))^{-1}$-recurrent.

Let $\tilde{\rho} > \mathfrak{a}^{-1}$, where $\mathfrak{a}$ is given as in (5.33), and apply Lemma 5.1. This states that, uniformly in $x \in \mathcal{L}_a h$ and $(\alpha, \varepsilon) \in \{(\alpha_j, \varepsilon_j)\}_{j \geq J}$, there exists a finite constant $D$ such that

$$
\begin{aligned}
D &\geq \sum_{n=1}^{\infty} \int_{\mathbb{S}} \tilde{\rho}^n \frac{r_\alpha^{(\varepsilon)}(y)}{r_\alpha^{(\varepsilon)}(x)} ((\hat{\mathfrak{P}}_\alpha^{(\varepsilon)})^{n-1} \underline{\hat{P}}_\alpha^{(\varepsilon)})(x, dy) \\
&= \sum_{n=1}^{\infty} \int_{\mathbb{S}} \tilde{\rho}^n (\theta^{(\varepsilon)}(\alpha))^n (\mathfrak{Q}_\alpha^{(\varepsilon)})^{n-1} Q_\alpha^{(\varepsilon)}(x, dy) \mathbf{1}_{\mathcal{L}_a h}(y) \\
&= \mathbb{E}_x^Q [(\tilde{\rho} \theta^{(\varepsilon)}(\alpha))^{\mathfrak{T}}],
\end{aligned}
$$
(5.64)

where $\mathbb{E}^Q[\cdot]$ denotes the expectation under the shifted measure $Q_\alpha^{(\varepsilon)}$. Now the probability measure $Q_\alpha^{(\varepsilon)}$ satisfies the minorization $(\mathcal{M}_Q)$ introduced in Section 5.1. It follows that upon each return to the set $\mathcal{L}_a h$, regeneration occurs with probability $g_\alpha^{(\varepsilon)}(x)$. Let $\hat{\mathbb{E}}_x^Q[\cdot]$ denote expectation conditioned on the event that regeneration does *not* occur at the start of a cycle evolving from $\mathcal{L}_a h$. Then as a consequence of (5.64), we obtain

$$(1 - g_\alpha^{(\varepsilon)}(x)) \hat{\mathbb{E}}_x^Q [(\tilde{\rho} \theta^{(\varepsilon)}(\alpha))^{\mathfrak{T}}] \leq D \tag{5.65}$$

uniformly in $x \in \mathcal{L}_a h$ and $(\alpha, \varepsilon) \in \{(\alpha_j, \varepsilon_j)\}_{j \geq J}$.

We begin by obtaining a lower bound for the function $g_\alpha^{(\varepsilon)}$. By the definition of $g_\alpha^{(\varepsilon)}$ [given in Section 5.1 under $(\mathcal{M}_Q)$],

$$g_\alpha^{(\varepsilon)}(x) \geq \gamma_\alpha^{(\varepsilon)} \wedge \frac{1}{2} \qquad \forall x \in \mathcal{L}_a h,$$

where

$$\gamma_\alpha^{(\varepsilon)} := \frac{\delta_a}{\theta^{(\varepsilon)}(\alpha)} \Big( \sup_{x \in \mathcal{L}_a h} r_\alpha^{(\varepsilon)}(x) \Big)^{-1} \int_{\mathbb{S}} e^{\alpha f(y)} r_\alpha^{(\varepsilon)}(y) \nu_a(dy)$$

and we now assert that $\gamma_\alpha^{(\varepsilon)} \geq \gamma_0$ for some positive constant $\gamma_0$. To establish this assertion, first recall by (5.33) that $\liminf_{j \to \infty} \theta^{(\varepsilon)}(\alpha_j) = \mathfrak{a} \in (\theta(\mathfrak{r}), 1]$, and by passing to an appropriate subsequence, we obtain the convergence as a *limit* rather than as a lower limit; that is,

$$(5.66) \qquad \lim_{j \to \infty} \theta^{(\varepsilon_j)}(\alpha_j) = \mathfrak{a} \in (\theta(\mathfrak{r}), 1].$$

Next recall by Proposition 5.2 that, for all $(\alpha, \varepsilon) \in \{(\alpha_j, \varepsilon_j)\}$, $r_\alpha^{(\varepsilon)}$ is bounded from above on $\mathcal{L}_a h$. Therefore, there exists a finite constant $b$ such that

$$(5.67) \qquad \liminf_{j \to \infty} \gamma_{\alpha_j}^{(\varepsilon_j)} \geq \frac{\delta_a}{b} \liminf_{j \to \infty} \int_{\mathbb{S}} e^{\alpha_j f(y)} r_{\alpha_j}^{(\varepsilon_j)}(y) \nu_a(dy).$$

Moreover, from the series representation of $r_\alpha^{(\varepsilon)}$ [given in its definition in Section 5.1], we also obtain that

$$(5.68) \qquad \liminf_{j \to \infty} r_{\alpha_j}^{(\varepsilon_j)}(x) \geq \frac{\delta_a \mathbf{1}_{\mathcal{L}_a h}(x)}{\mathfrak{a}} \qquad \forall x \in \mathbb{S}.$$

Then applying Fatou's lemma to the last term on the right-hand side of (5.67) yields

$$(5.69) \qquad \liminf_{j \to \infty} \int_{\mathbb{S}} e^{\alpha_j f(y)} r_{\alpha_j}^{(\varepsilon_j)}(y) \nu_a(dy) \geq \frac{\delta_a}{\mathfrak{a}} \int_{\mathcal{L}_a h} e^{\mathfrak{r} f(y)} \nu_a(dy) > 0,$$

where the last inequality follows since $\mathrm{supp}\, \nu_a \subseteq \mathcal{L}_a h$; cf. Remark 2.1. Thus we conclude that—possibly after passing to a subsequence—$\{\gamma_{\alpha_j}^{(\varepsilon_j)}\}$ is bounded from below by some positive constant, which we call $\gamma_0$, and we may assume that this constant has been chosen sufficiently small such that $\gamma_0 \in (0, 1/2]$.

Now choose $t > 1$ sufficiently large such that

$$(5.70) \qquad D < \frac{1}{2(1 - \gamma_0)^t}.$$

Then choose $(\alpha, \varepsilon) \in \{(\alpha_j, \varepsilon_j)\}_{j \in J}$ and $\rho > 1$ sufficiently small such that

$$(5.71) \qquad \frac{\rho^t}{\mathfrak{a}} \Big( \frac{\theta^{(\varepsilon)}(\alpha)}{\mathfrak{a}} \Big)^{t-1} \leq \tilde{\rho}.$$

[This choice is possible due to (5.66) and the fact that we have chosen $\tilde{\rho} > \mathfrak{a}^{-1}$. The point $(\alpha, \varepsilon)$ and constant $\rho$ satisfying (5.71) will now be fixed for the remainder of the proof.] Since $g_\alpha^{(\varepsilon)}(x) \leq 1/2$, substituting (5.71) into (5.65) gives

$$(5.72) \qquad 2D \geq \hat{\mathbb{E}}_x^Q\left[\left(\frac{\rho\theta^{(\varepsilon)}(\alpha)}{\mathfrak{a}}\right)^{t\mathfrak{T}}\right] \geq \hat{\mathbb{E}}_x^Q\left[\left(\frac{\rho\theta^{(\varepsilon)}(\alpha)}{\mathfrak{a}}\right)^{\mathfrak{T}}\right]^t \qquad \forall x \in \mathcal{L}_a h,$$

where the last step was obtained by Hölder's inequality. Then (5.70) and (5.72) yield

$$(5.73) \qquad D_0 := \sup_{x \in \mathcal{L}_a h} \hat{\mathbb{E}}_x^Q\left[\left(\frac{\rho\theta^{(\varepsilon)}(\alpha)}{\mathfrak{a}}\right)^{\mathfrak{T}}\right] < \frac{1}{1-\gamma_0}.$$

Our next objective is to show that if

$$\zeta := \frac{\rho\theta^{(\varepsilon)}(\alpha)}{\mathfrak{a}},$$

then (5.73) implies that $\mathbb{E}_x^Q[\zeta^\tau] < \infty$, all $x \in \mathcal{L}_a h$. To this end, let $N$ denote the random number of returns to the set $\mathcal{L}_a h$ which occur before regeneration actually takes place; that is to say, regeneration occurs directly following the $N$th visit to $\mathcal{L}_a h$. Then

$$(5.74) \qquad \sup_{x \in \mathcal{L}_a h} \mathbb{E}_x^Q[\zeta^\tau | N] \leq \left(\sup_{x \in \mathcal{L}_a h} \hat{\mathbb{E}}_x^Q[\zeta^{\mathfrak{T}}]\right)^N = D_0^N.$$

Now the probability that regeneration occurs upon any given return to $\mathcal{L}_a h$ is given by $g_\alpha(x) \geq \gamma_0$. Thus (5.73) and (5.74) yield

$$(5.75) \qquad \sup_{x \in \mathcal{L}_a h} \mathbb{E}_x^Q[\zeta^\tau] \leq \sum_{n=1}^{\infty} (1-\gamma_0)^{n-1} D_0^n < \infty.$$

Moreover, since $\theta^{(\varepsilon)}(\alpha) \searrow \mathfrak{a}$ implies that $\zeta > 1$, we have

$$(5.76) \qquad \mathbb{E}_x^Q\left[\sum_{n=1}^{\tau-1} \zeta^n\right] \leq \mathbb{E}_x^Q\left[\zeta^\tau\left(\frac{1}{\zeta} + \frac{1}{\zeta^2} + \cdots\right)\right] = \text{const.} \cdot \mathbb{E}_x^Q[\zeta^\tau]$$
$$\forall x \in \mathcal{L}_a h.$$

Therefore, since $\mathfrak{a} \leq 1$,

$$(5.77) \qquad \sup_{x \in \mathcal{L}_a h} \mathbb{E}_x^Q\left[\sum_{n=1}^{\tau-1} (\rho\theta^{(\varepsilon)}(\alpha))^n\right] \leq \sup_{x \in \mathcal{L}_a h} \mathbb{E}_x^Q\left[\sum_{n=1}^{\tau-1} \zeta^n\right] < \infty.$$

Finally observe that

$$\mathbb{E}_x^Q[(\rho\theta^{(\varepsilon)}(\alpha))^n \mathbf{1}_{\{\tau > n\}}] = (\rho\theta^{(\varepsilon)}(\alpha))^n (Q_\alpha^{(\varepsilon)} - g_\alpha^{(\varepsilon)} \otimes \mu_\alpha^{(\varepsilon)})^n (x, \mathbb{S}).$$

More explicitly, under the minorization ($\mathcal{M}_Q$), $1 - g_\alpha^{(\varepsilon)}$ describes the probability that regeneration does not occur at any given time step, and $(Q_\alpha^{(\varepsilon)} -$ $g_\alpha^{(\varepsilon)} \otimes \mu_\alpha^{(\varepsilon)})/(1 - g_\alpha^{(\varepsilon)})$ describes the transition kernel in the event that regeneration does not occur. It then follows from the series representation of a regeneration cycle [as utilized, e.g., in Nummelin (1978) or Ney and Nummelin (1987a), Lemma 4.1] that

$$
\begin{aligned}
(5.78) \quad \mathbb{E}_x^Q\left[\sum_{n=1}^{\tau-1} (\rho\theta^{(\varepsilon)}(\alpha))^n\right] &= \sum_{n=1}^{\infty} (\rho\theta^{(\varepsilon)}(\alpha))^n (Q_\alpha^{(\varepsilon)} - g_\alpha^{(\varepsilon)} \otimes \mu_\alpha^{(\varepsilon)})^n(x,\mathbb{S}) \\
&\geq \frac{1}{r_\alpha^{(\varepsilon)}(x)} \sum_{n=1}^{\infty} \rho^n((\hat{P}_\alpha - \delta_a \mathbf{1}_{\mathcal{L}_a h} \otimes \hat{\nu}_a^{(\alpha)})^n r_\alpha^{(\varepsilon)})(x) \\
&= \mathbb{E}_x\left[\sum_{n=1}^{\tau-1} \frac{r_\alpha^{(\varepsilon)}(X_n)}{r_\alpha^{(\varepsilon)}(X_0)} \rho^n e^{\alpha S_n}\right].
\end{aligned}
$$

[The inequality comes from the fact that, in the definition of $g_\alpha^{(\varepsilon)}$, we have truncated this quantity at the level $1/2$. Also, inequality results since we have substituted $\hat{P}_\alpha$ for $\hat{P}_\alpha^{(\varepsilon)}$.] Hence

$$
(5.79) \qquad \sup_{x\in\mathcal{L}_a h} \mathbb{E}_x\left[\sum_{n=1}^{\tau-1} \frac{r_\alpha^{(\varepsilon)}(X_n)}{r_\alpha^{(\varepsilon)}(X_0)} \rho^n e^{\alpha S_n}\right] < \infty.
$$

The eigenvectors can effectively be removed from this last expression, since $(\alpha,\varepsilon)$ may now be fixed, and then $r_\alpha^{(\varepsilon)}$ is uniformly bounded from below by a positive constant, while $r_\alpha^{(\varepsilon)}$ is uniformly bounded from above on $\mathcal{L}_a h$. Moreover, $\text{supp}\,\nu_a \subseteq \mathcal{L}_b f$ (cf. Remark 2.1), so $X_0 \sim \nu_a \Longrightarrow f(X_0) \leq b < \infty$. Hence $S_n' \leq S_n + b$, all $n \geq 0$. Therefore it follows from (5.79) that

$$
(5.80) \qquad \mathbb{E}_{\nu_a}\left[\sum_{n=0}^{\tau-1} \rho^n e^{\alpha S_n'}\right] < \infty
$$

as required.

If $\hat{P}_\alpha^{(\varepsilon)}$ is $(\theta^{(\varepsilon)}(\alpha))^{-1}$-transient, then we can normalize the measure $Q_\alpha^{(\varepsilon)}$, multiplying it by an appropriate constant $\lambda > 1$, so that $\lambda Q_\alpha^{(\varepsilon)}$ is a probability measure. Let

$$
R_\alpha^{(\varepsilon)}(x,E) = \lambda Q_\alpha^{(\varepsilon)}(x,E) \quad \text{and} \quad \mathfrak{R}_\alpha^{(\varepsilon)}(x,E) = \lambda \mathfrak{Q}_\alpha^{(\varepsilon)}(x,E) \qquad \forall x,E.
$$

Since $\lambda > 1$, the measure $R_\alpha^{(\varepsilon)}$ satisfies the same minorization as $Q_\alpha^{(\varepsilon)}$, namely,

$$
(\mathcal{M}_R) \qquad g_\alpha^{(\varepsilon)}(x)\mu_\alpha^{(\varepsilon)}(E) \leq R_\alpha^{(\varepsilon)}(x,E) \qquad \forall x\in\mathbb{S},\ E\in\mathcal{S};
$$

cf. $(\mathcal{M}_R)$ in Section 5.1. Hence the regularity properties that have just been developed for $g_\alpha^{(\varepsilon)}$ may also be applied to the present case without any modification.

To obtain (5.63), choose $\tilde{\rho} > \lambda/\mathfrak{a}$ and apply Lemma 5.1, just as in (5.64), (5.65), to obtain

$$(5.81) \qquad D \geq \mathbb{E}_x^R\left[\left(\frac{\tilde{\rho}\theta^{(\varepsilon)}(\alpha)}{\lambda}\right)^{\mathfrak{T}}\right] \geq (1 - g_\alpha^{(\varepsilon)}(x))\hat{\mathbb{E}}_x^R\left[\left(\frac{\tilde{\rho}\theta^{(\varepsilon)}(\alpha)}{\lambda}\right)^{\mathfrak{T}}\right],$$

uniformly in $x \in \mathcal{L}_a h$, where $\hat{\mathbb{E}}_x^R[\cdot]$ denotes expectation under $R$ conditioned on the event that regeneration does not occur at the start of the regeneration cycle.

Choose $t \geq 1$ sufficiently large so that (5.70) holds, and then choose $(\alpha, \varepsilon)$ and $\rho > \lambda$ sufficiently small such that

$$(5.82) \qquad \frac{\rho^t}{\mathfrak{a}}\left(\frac{\theta^{(\varepsilon)}(\alpha)}{\lambda\mathfrak{a}}\right)^{t-1} \leq \tilde{\rho}.$$

[As $(\alpha, \varepsilon) \searrow (\mathfrak{r}, 0)$ and $\rho \searrow \lambda$, the left-hand side of (5.82) converges to $\lambda/\mathfrak{a} < \tilde{\rho}$, where the last equality follows from the original choice of $\tilde{\rho}$. Thus, it is always possible to find appropriate elements $(\alpha, \varepsilon)$ and $\rho$ satisfying (5.82).] Using (5.82) in place of (5.71), we then obtain (5.72), but with $\hat{\mathbb{E}}_x^Q[\cdot]$ replaced everywhere with $\hat{\mathbb{E}}_x^R[\cdot]$ and "$\theta^{(\varepsilon)}(\alpha)$" replaced everywhere with "$\theta^{(\varepsilon)}(\alpha)/\lambda$." The remainder of the proof may now be repeated without change [except that "$\theta^{(\varepsilon)}(\alpha)$" is replaced everywhere with "$\theta^{(\varepsilon)}(\alpha)/\lambda$"] to obtain (5.77) or, more precisely,

$$(5.83) \qquad \sup_{x \in \mathcal{L}_a h} \mathbb{E}_x^R\left[\sum_{n=1}^{\tau-1}\left(\frac{\rho\theta^{(\varepsilon)}(\alpha)}{\lambda}\right)^n\right] < \infty.$$

To complete the proof, note that

$$
\begin{aligned}
(5.84) \qquad \mathbb{E}_x^R\left[\sum_{n=1}^{\tau-1}\left(\frac{\rho\theta^{(\varepsilon)}(\alpha)}{\lambda}\right)^n\right] &= \sum_{n=1}^{\infty}\left(\frac{\rho\theta^{(\varepsilon)}(\alpha)}{\lambda}\right)^n (R_\alpha^{(\varepsilon)} - g_\alpha^{(\varepsilon)} \otimes \mu_\alpha^{(\varepsilon)})^n(x, \mathbb{S}) \\
&\geq \sum_{n=1}^{\infty}(\rho\theta^{(\varepsilon)}(\alpha))^n (Q_\alpha^{(\varepsilon)} - g_\alpha^{(\varepsilon)} \otimes \mu_\alpha^{(\varepsilon)})^n(x, \mathbb{S}) \\
&\geq \mathbb{E}_x\left[\sum_{n=1}^{\tau-1} \frac{r_\alpha^{(\varepsilon)}(X_n)}{r_\alpha^{(\varepsilon)}(X_0)}\rho^n e^{\alpha S_n}\right].
\end{aligned}
$$

An inequality arises in the second equation when, on the right-hand side, we replace $\lambda^{-1}g_\alpha^{(\varepsilon)}(x)$ with the larger quantity $g_\alpha^{(\varepsilon)}(x)$. The final inequality is then obtained as in (5.78). Consequently we conclude that (5.79) holds and hence also (5.80). $\quad\square$

PROOF OF THEOREM 4.2. Set

$$S_n' = S_n + f(X_0) \qquad \forall n \geq 0,$$

and

$$\check{S} \stackrel{d}{=} S_{T_i - 1} - S_{T_{i-1} - 1} \qquad \forall i \geq 1.$$

Then by Lemma 5.2,

$$(5.85) \qquad \mathbb{E}[\check{A}^\alpha] := \mathbb{E}[e^{\alpha\check{S}}] < \infty \qquad \text{for some } \alpha > \mathfrak{r},$$

which establishes the first assertion of the theorem.

We now turn to the remaining estimates stated in (4.16). Assume that regeneration occurs at time zero and let $\tau$ denote the subsequent regeneration time. Let

$$\check{M}^* = \sup\{|B_0| + A_0|B_1| + \cdots + (A_0 \cdots A_{j-1})|B_j| : 0 \leq j < \tau\}$$

and observe that

$$(5.86) \qquad \mathbb{E}[|\check{B}|^\alpha] \leq \mathbb{E}[(\check{M}^*)^\alpha] \leq \mathbb{E}_{\nu_a}\left[\left(|B_0| + \sum_{n=1}^\infty e^{S'_{n-1}}|B_n|\mathbf{1}_{\{\tau > n\}}\right)^\alpha\right].$$

Now if $\alpha \geq 1$, then by Minkowski's inequality and the independence of $\{B_n\}$ from $(\tau, \{A_n\})$,

$$\mathbb{E}_{\nu_a}\left[\left(|B_0| + \sum_{n=1}^\infty e^{S'_{n-1}}|B_n|\mathbf{1}_{\{\tau > n\}}\right)^\alpha\right]^{1/\alpha}$$

$$(5.87) \qquad \leq \mathbb{E}[|B|^\alpha]^{1/\alpha} + \sum_{n=1}^\infty \mathbb{E}[|B|^\alpha]^{1/\alpha} \mathbb{E}_{\nu_a}[e^{\alpha S'_{n-1}}\mathbf{1}_{\{\tau > n\}}]^{1/\alpha}$$

$$\leq (1 + \mathbb{E}[|B|^\alpha])\left(1 + \sum_{n=0}^\infty \mathbb{E}_{\nu_a}[e^{\alpha S'_n}\mathbf{1}_{\{\tau > n+1\}}]^{1/\alpha}\right) < \infty,$$

where the final inequality holds for sufficiently small $\alpha > \mathfrak{r}$, since (H$_2$) yields $\mathbb{E}[|B|^\alpha] < \infty$, while Lemma 5.2 yields

$$(5.88) \qquad \mathbb{E}_{\nu_a}[e^{\alpha S'_n}\mathbf{1}_{\{\tau > n\}}] \leq K\rho^{-n}, \qquad n = 0, 1, \ldots,$$

for some constants $K < \infty$ and $\rho > 1$.

On the other hand, if $\alpha < 1$, then in place of Minkowski's inequality we may apply the deterministic inequality

$$(5.89) \qquad (y + z)^\gamma \leq y^\gamma + z^\gamma \qquad \forall y \geq 0, \ z \geq 0 \text{ and } 0 \leq \gamma \leq 1,$$

and the result follows in the same way as before.

Finally consider (4.17). Assume, to the contrary, that there exists a $\varphi$-positive set $F \in \mathcal{S}$ and a finite constant $b$ such that

$$(5.90) \qquad f(x) \geq -b \quad \text{and} \quad \mathbb{E}_x[\check{M}_0^\alpha] = \infty \qquad \forall x \in F.$$

Note

$$\check{M}_0^* := \sup\{|B_0| + A_0|B_1| + \cdots + (A_0 \cdots A_{j-1})|B_j| : 0 \leq j < T_0\}$$
$$\geq |B_0| + A_0 \sup\{B_1 + A_1 B_2 + \cdots + (A_1 \cdots A_{j-1})B_j : 1 \leq j < T_0\}$$
$$:= |B_0| + A_0 \check{M}_0.$$

If $\alpha > 1$, then it follows by Minkowski's inequality and (5.90) that

$$(5.91) \quad \mathbb{E}_x[(\check{M}_0^*)^\alpha]^{1/\alpha} \geq -\mathbb{E}[|B|^\alpha]^{1/\alpha} + e^{-b/\alpha}\mathbb{E}_x[\check{M}_0^\alpha]^{1/\alpha} \qquad \forall x \in F.$$

Then by (5.90),

$$(5.92) \qquad \mathbb{E}_x[(\check{M}_0^*)^\alpha] = \infty \qquad \forall x \in F.$$

If $\alpha \leq 1$, then we may apply (5.89) in place of Minkowski's inequality, and the conclusion still holds.

Now suppose that regeneration occurs at time zero, and set

$$\mathcal{T} = \inf\{n \in \mathbb{Z}_+ : X_n \in F\} \wedge \tau.$$

Since $\varphi(F) > 0$, $\{X_n\}$ will visit the set $F$ over the given regeneration cycle with positive probability. Since $\mathbb{E}[\check{M}^\alpha] < \infty$ implies $\nu_a(F) = 0$, this means that $\mathbf{P}\{\mathcal{T} < \tau\} > 0$. Let

$$\check{N}^* = \sup\{|B_{\mathcal{T}}| + A_{\mathcal{T}}|B_{\mathcal{T}+1}| + \cdots + (A_{\mathcal{T}} \cdots A_{j-1})|B_j| : \mathcal{T} \leq j < \tau\}$$

and observe by definition that

$$(5.93) \qquad \check{M}^* \geq (A_1 \cdots A_{\mathcal{T}-1})\check{N}^*.$$

Now conditional on the event $\mathcal{E} := \{\mathcal{T} < \tau\}$,

$$(5.94) \quad \mathbb{E}_\mathcal{E}[\mathbb{E}_\mathcal{E}[(A_1 \cdots A_{\mathcal{T}-1}\check{N}^*)^\alpha | \mathcal{T}, A_1, \ldots, A_{\mathcal{T}-1}, X_\mathcal{T}]]$$
$$= \mathbb{E}_\mathcal{E}[(A_1 \cdots A_{\mathcal{T}-1})^\alpha \mathbb{E}_{X_\mathcal{T}}[(\check{N}^*)^\alpha]].$$

Moreover, the left-hand side of (5.94) is finite, since (5.93) holds and we have shown that $\mathbb{E}[(\check{M}^*)^\alpha] < \infty$. However, on the right-hand side of (5.94), the term $\mathbb{E}_{X_\mathcal{T}}[(\check{N}^*)^\alpha]$ is infinite due to (5.92). (In particular, note that the definitions of $\check{N}^*$ and $\check{M}_0^*$ are the same, except that the quantities in the definition of $\check{N}^*$ are conditional on an initial state $X_\mathcal{T} \in F$, while the quantities in the definition of $\check{M}_0^*$ are conditional on an initial state $X_0 = x \in \mathbb{S}$.) Consequently, we conclude

$$\mathbb{E}_\mathcal{E}[(A_1 \cdots A_{\mathcal{T}-1})^\alpha] = 0.$$

But since $\mathcal{T} \leq \tau < \infty$ a.s. and $A_n > 0$ for all $n$, this is impossible. Hence $\mathbb{E}_x[\check{M}_0^\alpha] < \infty$ for $\varphi$-a.a. $x$, and then $\mathbb{E}_x[\check{B}_0^\alpha] \leq \mathbb{E}_x[\check{M}_0^\alpha] < \infty$, $\varphi$-a.a. $x$. Setting $B_n \equiv 1$ and observing that in this case $\check{B}_0 \geq \check{A}_0$, we also obtain $\mathbb{E}_x[\check{A}_0^\alpha] < \infty$, $\varphi$-a.a. $x$. $\quad\square$

**6. Some extensions.**

6.1. *A general Markovian model.* In this section, we consider an extension which allows the sequence $\{B_n\}$ to be Markov dependent, and for there to be dependence between $\{A_n\}$ and $\{B_n\}$. Now assume

$$(6.1) \qquad \log A_n = f(X_n) \quad \text{and} \quad B_n = G(X_n),$$

where $f : \mathbb{S} \to \mathbb{R}$ and $G : \mathbb{S} \to \mathbb{R}$. For simplicity, suppose for the moment that $|B_n| \geq 1$ for all $n$, and set $g(x) = \log G(x)$. Let

$$\tilde{S}_n = \log A_1 + \cdots + \log A_{n-1} + \log |B_n|,$$

and let $\tilde{\boldsymbol{\Lambda}}$ be defined as $\boldsymbol{\Lambda}$, but with $\tilde{S}_n$ in place of $S_n$. In addition to (H$_2$), assume that $\tilde{\boldsymbol{\Lambda}}(\boldsymbol{\alpha}) < \infty$ for some $\boldsymbol{\alpha} = (\alpha, \beta)$, where $\alpha > \mathfrak{r}$ and $\beta > 0$.

An approach to this more general problem is to introduce the kernel

$$(6.2) \qquad \hat{K}_\alpha(x, E) := \int_E e^{\alpha F(x,y)} P(x, dy) \qquad \forall x \in \mathbb{S}, \ E \in \mathcal{S},$$

where $F(x, y) = f(x) + (g(y) - g(x))$. Note from this definition that

$$(6.3) \qquad \hat{K}_\alpha^n(x, \mathbb{S}) = \mathbb{E}_x[|B_0|^{-\alpha}(A_0 \cdots A_{n-1})^\alpha |B_n|^\alpha],$$

that is, $\hat{K}_\alpha$ now plays the role of $\hat{P}_\alpha$, and similarly for the kernel $\mathfrak{K}(x, E) := \int_{E \cap (\mathcal{L}_a h)^c} \hat{K}_\alpha(x, dy)$, $\mathfrak{K}^{(\varepsilon)}$, and so on. We note that the random quantity $|B_0|$ may be bounded under an appropriate choice of supp $\nu_a$, as in Remark 2.1. In this way we obtain as a direct extension of Lemma 5.2 that

$$(6.4) \qquad \mathbb{E}_{\nu_a}\left[\left(|B_0| + \sum_{n=1}^{\tau-1} \rho^{n-1} e^{S_{n-1}} |B_n|\right)^\alpha\right] < \infty,$$

which is the critical estimate needed to handle the present generalization. Finally, the assumption $|B_n| \geq 1$ may be dropped by replacing $|B_n|$ with $|B_n| + 1$ in the relevant parts of the proofs.

In practice, dependence arises in the sequence $\{B_n\}$ if, for example, the premiums of the insurance company are determined by a bonus system, in which case they depend on the observed claims in the previous time intervals. In this setting, $\{B_n\}$ may be modeled as a function of a Markov chain on $\mathbb{R}$, say. See Lemaire (1995) for general background on bonus systems, or Bonsdorff (2005) for some recent developments.

6.2. *Generalizations of* $(\mathfrak{M})$. Finally we discuss certain generalizations of the minorization condition $(\mathfrak{M})$. A simple extension—still general enough to handle, for example, the AR($p$) model when $p > 1$—is to replace the Markov chain $\{X_n\}$ with

$$(6.5) \qquad \mathcal{X}_n := (X_{(n-1)k+1}, \ldots, X_{nk}), \qquad n = 0, 1, \ldots,$$

for any positive integer $k$, and to introduce the following minorization condition:

$(\mathcal{M}_1)$   $\delta_a \mathbf{1}_{\mathcal{L}_a h}(x)\nu_a(E) \leq \mathbb{P}\{\mathcal{X}_1 \in E | \mathcal{X}_0 = x\}$   $\forall x \in \mathbb{S}^k,\ E \in \mathcal{S}^k.$

Set $\mathcal{A}_n = A_{(n-1)k+1} \cdots A_{nk}$ and assume that the analogs of (H$_1$)–(H$_3$) hold. Then a repetition of the argument given in Lemma 5.2 yields

$$(6.6) \qquad \mathbb{E}[\mathcal{A}_0^\alpha + \rho \mathcal{A}_0^\alpha \mathcal{A}_1^\alpha + \cdots + \rho^{\tau-1}\mathcal{A}_0^\alpha \cdots \mathcal{A}_{\tau-1}^\alpha] < \infty$$

for some $\alpha > \mathfrak{r}$ and $\rho > 1$,

where $\tau$ denotes a typical regeneration time of the chain $\{\mathcal{X}_n\}$. From (6.6) we immediately obtain the estimate for $\check{A}$ given in Theorem 4.2 [although not the estimates for $\check{B}$ or $\check{M}$, since (6.6) omits terms from the sequence on the left-hand side of (5.63)]. To obtain the remaining estimates of that theorem, we could, for example, introduce the random variable

$$(6.7) \qquad \mathcal{B}_n = B_{(n-1)k+1} + \cdots + A_{(n-1)k+1} \cdots A_{nk-1} B_{nk}$$

and follow the approach just outlined in the previous section. As an alternative and more direct approach, we could work explicitly with the original sequence $\{B_n\}$ and with $(\mathcal{M}_1)$, but introduce a strengthening of (H$_2$). Namely set $\log \bar{A}_n = \max\{\log A_n, 0\}$ and $\bar{\mathcal{A}}_n = \bar{A}_{(n-1)k+1} \cdots \bar{A}_{nk}$, and assume the further requirement:

(H$_2'$)   Hypothesis (H$_2$) holds with $\bar{S}_n := \log \bar{A}_1 + \cdots + \log \bar{A}_n$ in place of $S_n$.

[A precise statement of this condition would also involve a slight change in the definition of $h$; see the discussion below, where we explicitly verify this condition in the AR($p$) case.] A slight variant on our argument then yields

$$(6.8) \qquad \mathbb{E}[\bar{\mathcal{A}}_0^\alpha + \rho \mathcal{A}_0^\alpha \bar{\mathcal{A}}_1^\alpha + \cdots + \rho^{\tau-1}\mathcal{A}_0^\alpha \cdots \mathcal{A}_{\tau-2}^\alpha \bar{\mathcal{A}}_{\tau-1}^\alpha] < \infty$$

for some $\alpha > \mathfrak{r}$ and $\rho > 1$,

which is a mild extension of (6.6). If $\{B_n\}$ is i.i.d. and independent of $\{A_n\}$, then we consequently obtain the required estimates in Theorem 4.2 for $\check{B}$ and $\check{M}$ and hence the limit results stated in our two main theorems. In summary, if we assume the stronger condition (H$_2'$), then our main results hold under the alternative minorization condition $(\mathcal{M}_1)$.

For example, if $\{X_n\}$ is an AR($p$) process with $p > 1$, that is,

$$(6.9) \qquad X_n = \sum_{i=1}^p a_i X_{n-i} + \zeta_n, \qquad n = 1, 2, \ldots,$$

where $\{\zeta_n\}$ is an i.i.d. Gaussian sequence and $X_0 = x \in \mathbb{R}$, and if

$f(x) = x - \mu$   for some $\mu > 0$,

then ($\mathcal{M}_0$) will *not* hold with $k = 1$. However, the minorization ($\mathcal{M}_1$) and associated regularity conditions can be verified explicitly.

To demonstrate that our conditions actually hold in this case, we begin by expressing (6.9) in matrix form. Namely, for each $n \in \mathbb{Z}_+$ set

$$\mathbf{X}_n = \begin{pmatrix} X_n \\ \vdots \\ X_{n-p+1} \end{pmatrix}, \qquad \boldsymbol{\zeta}_n = \begin{pmatrix} \zeta_n \\ 0 \\ \vdots \\ 0 \end{pmatrix} \quad \text{and} \quad \mathfrak{A} = \begin{pmatrix} a_1 & \cdots & \cdots & a_p \\ 1 & & & 0 \\ & \ddots & & \vdots \\ 0 & & 1 & 0 \end{pmatrix}$$

(where $X_0, \ldots, X_{-p+1}$ are taken to be arbitrary deterministic values). It follows from these definitions that

(6.10)            $\mathbf{X}_n = \mathfrak{A}\mathbf{X}_{n-1} + \boldsymbol{\zeta}_n, \qquad n = 0, 1, \ldots.$

Set $\mathcal{X}_0 = \mathbf{X}_0, \mathcal{X}_1 = \mathbf{X}_p, \mathcal{X}_2 = \mathbf{X}_{2p}, \ldots.$

From (6.10) we obtain

(6.11)            $\mathbf{X}_n = \mathfrak{A}^n \mathbf{X}_0 + \mathfrak{A}^{n-1} \boldsymbol{\zeta}_1 + \cdots + \boldsymbol{\zeta}_n, \qquad n = 0, 1, \ldots.$

Hence, in particular,

(6.12)            $\mathbf{X}_p = \mathfrak{A}^p \mathbf{X}_0 + \mathbf{W}_p$

for some random vector $\mathbf{W}_p \sim \text{Normal}(0, S)$, where $S$ is a covariance matrix which is easily seen to have rank $p$.

From (6.12) it follows that for any $a > 0$, there exists a finite constant $b$ such that

(6.13)            $\mathbf{X}_0 \in \mathcal{B}_a(0) \quad \Longrightarrow \quad \mathfrak{A}^p \mathbf{X}_0 \in \mathcal{B}_b(0),$

where $\mathcal{B}_r(x)$ denotes a ball of radius $r$ about $x$. Let $\Phi_x$ denote the Normal$(x, S)$ density function, and let

$$\underline{\Phi}_b(y) = \inf\{\Phi_x(y) : x \in \mathcal{B}_b(0)\} \qquad \forall y \in \mathbb{R}^p.$$

Note that $\underline{\Phi}_b$ is positive everywhere since, for any fixed $y$, $\Phi_x(y)$ is continuous as a function of $x$ and hence achieves its minimum on the compact set $\mathcal{B}_b(0)$. Then by (6.13),

(6.14)            $\displaystyle\int_E \underline{\Phi}_b(y) \, dy \le P(x, E) \qquad \forall x \in \mathcal{B}_a(0), \ E \in \mathcal{S}^p.$

Consequently ($\mathcal{M}_1$) holds with $\nu_a(dy) = c\underline{\Phi}_b(y) \, dy$ and $c \in (0, \infty)$ a normalizing constant.

The verification of (H$_1$) is likewise straightforward. For example, the rate function $\mathfrak{r}$ appearing in (H$_1$) is just the nonzero point at which $\Lambda(\mathfrak{r}) = 0$, where $\Lambda$ is the Gärtner–Ellis limit for the process $\{S_n\}$, and this limit may be computed, as in Example 3.2, by observing that $S_n$ is normally distributed for all $n$, so it is sufficient to calculate the limiting values of its normalized mean and variance, which can be shown to converge as $n \to \infty$. In this way we obtain, more explicitly, that $\Lambda(\alpha) = -\alpha m + \sigma^2 \alpha^2 / 2$ for certain positive constants $m$ and $\sigma$ (the cumulant generating function for an appropriate normal distribution); thus, in particular, $\Lambda'(0) < 0$.

For (H$_3$), first recall that under an appropriate linear transformation, $\mathcal{L}$, we have that $\mathbf{W}_k \mapsto \tilde{\mathbf{W}}_k \sim \text{Normal}(0, I)$, where $I$ is the identity matrix. Let $\tilde{\Phi}$ denote the Normal$(0, I)$ density function, and fix $\tilde{a} > 0$. Then there exists a finite constant $b$ such that if $\tilde{x} \in \partial \mathcal{B}_b(0)$, then

$$\tilde{\Phi}_z(y) \le \tilde{\Phi}_{\tilde{x}}(y) \qquad \forall z \in \mathcal{B}_{\tilde{a}}(0), \ y \in \mathfrak{C}_\gamma(\tilde{x}), \tag{6.15}$$

where

$$\mathfrak{C}_\gamma(\tilde{x}) := \{\omega z : \omega \ge 1, z \in \mathcal{B}_\gamma(\tilde{x})\}$$

is a $\gamma$-cone eminating from $\tilde{x}$ and $\gamma > 0$ is a sufficiently small constant. Thus there exists a finite collection of points $\{\tilde{x}_1, \ldots, \tilde{x}_l\} \in \partial \mathcal{B}_b(0)$ such that

$$\tilde{P}(z, E) \le \sum_{i=1}^{l} \tilde{P}(\tilde{x}_i, E) \qquad \forall z \in \mathcal{B}_{\tilde{a}}(0), \ \forall E,$$

where $\tilde{P}$ denotes the transition kernel of $\{\mathcal{X}_n\}$ under the transformation $\mathcal{L}$. Under the inverse transformation, $\mathcal{L}^{-1}$, we then obtain that for any given $a > 0$, there exist points $\{x_1, \ldots, x_l\}$ lying on the boundary of some ellipse such that

$$P(z, E) \le \sum_{i=1}^{l} P(x_i, E) \qquad \forall z \in \mathcal{B}_a(0), \ \forall E, \tag{6.16}$$

and by a slight extension, we may replace $x_i$ with $\mathcal{B}_\varepsilon(x_i)$ (for some $\varepsilon > 0$) on the right-hand side.

The remaining condition to be verified is that (H$_2$) holds with

$$\bar{S}_n := \log \bar{A}_1 + \cdots + \log \bar{A}_n$$

in place of $S_n$, where $\log \bar{A}_n := \max\{\log A_n, 0\}$. Since $\log A_n = f(X_n) = X_n - \mu$ for some $\mu > 0$, we clearly have

$$\bar{S}_n \le |X_1| + \cdots + |X_n|.$$

To characterize the quantity on the right-hand side, note as a consequence of (6.11) that

$$|X_n| \le \tilde{D}(\lambda^n |X_0| + \lambda^{n-1} |\zeta_1| + \cdots + |\zeta_n|) \qquad \forall n, \tag{6.17}$$

where $\lambda$ is the spectral radius of $\mathfrak{A}$ and $\tilde{D}$ is a finite constant. Now under the standard conditions needed to ensure stationarity of the Markov chain in (6.9) [cf. Brockwell and Davis (1991)], the spectral radius of $\mathfrak{A}$ is less than one. Hence

$$(6.18) \qquad \bar{S}_n \le D(|X_0| + |\zeta_1| + \cdots + |\zeta_n|)$$

for some finite constant $D$. Since the sequence $\{\zeta_n\}$ is i.i.d., it follows that

$$(6.19) \qquad \mathbb{E}[e^{(\alpha/D)|\zeta_1|}] < \infty \quad \Longrightarrow \quad \limsup_{n\to\infty} \frac{1}{n} \log \mathbb{E}[e^{\alpha \bar{S}_n}] < \infty.$$

As $\zeta_1$ has a standard Gaussian distribution, we conclude

$$(6.20) \qquad \limsup_{n\to\infty} \frac{1}{n} \log \mathbb{E}[e^{\alpha \bar{S}_n}] < \infty \qquad \forall \alpha.$$

Now in a precise statement of $(\mathrm{H}_2)$, we would also need to specify a function $h\colon \mathbb{S}^p \to \mathbb{R}$ which corresponds to the Markov chain $\{\mathcal{X}_n\}$ rather than to $\{X_n\}$. What is actually needed is that

$$(6.21) \qquad \limsup_{n\to\infty} \frac{1}{n} \log \mathbb{E}[e^{\alpha \bar{S}_{pn} + \beta S_n^{(h)}}] < \infty, \qquad \text{some } \alpha > \mathfrak{r} \text{ and } \beta > 0,$$

where $S_n^{(h)} = h(\mathcal{X}_1) + \cdots + h(\mathcal{X}_n)$. In the above discussion, we have chosen

$$h(x) = \|x\| \le |x_1| + \cdots + |x_p|, \qquad \text{where } x = (x_1, \ldots, x_p).$$

But then (6.21) follows from (6.20) [and its proof, since we may replace $\bar{S}_n$ with $\sum_{i=1}^n |X_n|$ in the deduction following (6.18)].

From these considerations, we conclude that if $\{X_n\}$ is an AR($p$) process with $p > 1$, then the results of this paper still apply; in particular, the examples in Section 3 can all be considered in this more general setting.

Finally, in the general ARMA($p,q$) case, where

$$(6.22) \qquad X_n = \sum_{i=1}^p a_i X_{n-i} + \sum_{j=0}^q b_j \zeta_{n-j}, \qquad n = 1, 2, \ldots,$$

it may easily be shown [cf. Meyn and Tweedie (1993), page 28] that

$$(6.23) \qquad X_n = b_0 Y_n + \cdots + b_q Y_{n-q},$$

where $\{Y_n\}$ is the corresponding AR($p$) process obtained by setting each $b_j$ in (6.22) to zero. Thus letting

$$\mathcal{Y}_n = (Y_{n(p-1)+1}, \ldots, Y_{np}) \quad \text{and} \quad \mathcal{A}_n = A_{(n-1)p+1} \cdots A_{np},$$

we see that $\mathcal{A}_n = F(\mathcal{Y}_{n-1}, \mathcal{Y}_n)$ (where, if necessary, we set $c_j = 0$ to force $p \ge q + 1$).

In this and more general situations, it is useful to consider a further extension of $(\mathfrak{M})$, as follows. Let $\{\mathcal{A}_n\}$ again be defined as in the discussion prior to (6.6) (i.e., with index "$k$" in place of "$p$"), and let $\{\mathcal{B}_n\}$ be defined as in (6.7). Let $V_n = h(\mathcal{X}_{n-1}, \mathcal{X}_n)$, and assume that

$$\xi_n := (\mathcal{A}_n, \mathcal{B}_n, V_n) = F(\mathcal{X}_{n-1}, \mathcal{X}_n),$$

for some function $F : \mathbb{S}^{2k} \to \mathbb{R}^3$. Next introduce the following minorization condition:

$(\mathcal{M}_2)$     $\delta_a \mathbf{1}_{\mathcal{L}_a h}(x) \nu_a(E \times \Gamma) \leq \mathbb{P}\{(\mathcal{X}_1, \xi_1) \in E \times \Gamma | \mathcal{X}_0 = x\}$

$$\forall x \in \mathbb{S}^k, \ E \in \mathcal{S}^k, \ \Gamma \in \mathfrak{R}^3.$$

Using a result in Ney and Nummelin [(1986), page 4] in place of Lemma 4.1, we obtain under $(\mathcal{M}_2)$ that $\{(\mathcal{X}_n, \xi_n)\}$ exhibits a regeneration structure and, moreover, that its transform kernel exhibits a minorization of the form $(\hat{\mathcal{M}})$ [where $(\hat{\mathcal{M}})$ is given as in the proof of Proposition 5.1]. Then one can proceed as before, although a rigorous analysis requires a careful treatment and further moment conditions on the process $(\mathcal{A}_n, \mathcal{B}_n, V_n)$, which we do not pursue here.

**Acknowledgments.**   The author is particularly grateful to Harri Nyrhinen and Esa Nummelin for many fruitful discussions and for their hospitality during several visits to Univ. Helsinki. The author would also like to thank the referees for many helpful suggestions and improvements.

DEPARTMENT OF MATHEMATICAL SCIENCES
UNIVERSITY OF COPENHAGEN
UNIVERSITETSPARKEN 5
DK-2100 COPENHAGEN Ø
DENMARK
E-MAIL: collamore@math.ku.dk